\documentclass[11pt]{amsart}

\usepackage{amsmath,amssymb,amsthm,mathrsfs,eucal,textcomp}                                                        
\usepackage[utf8]{inputenc}
\usepackage[T1]{fontenc}
\usepackage{lmodern}
\usepackage{amsmath}
\usepackage{amsfonts}
\usepackage{amssymb}
\usepackage{caption}
\usepackage{mathrsfs}
\usepackage{amsthm}
\usepackage{tikz}
\usepackage{bm}
\usepackage{color}
\usepackage{rotating}
\usepackage{array}

\theoremstyle{plain}
\newtheorem{theorem}{Theorem}[section]
\newtheorem{lemma}[theorem]{Lemma}
\newtheorem{proposition}[theorem]{Proposition}
\theoremstyle{definition}
\newtheorem{definition}[theorem]{Definition}
\newtheorem{remark}[theorem]{Remark}

\def\R{\mathbb{R}}
\def\e{\varepsilon}
\def\W{\mathcal{W}}
\def\PP{\mathcal{P}_2(\R)}

\def\K{\mathcal{K}}
\newcommand{\Abs}{\mathrm{AC}}
\newcommand{\sign}{\mathrm{sign}}
\def\leb{\mathcal{L}}

\def\mutil{\widetilde{\mu}}
\newcommand*\mycirc[1]{
  \begin{tikzpicture}
    \node[draw,circle,inner sep=0.8pt] {#1};
  \end{tikzpicture}}

\usepackage[scale=0.775]{geometry}

\title[Equivalence of gradient flows and entropy solutions]{Equivalence of gradient flows and entropy solutions for singular nonlocal interaction equations in 1D}

\author{G. A. Bonaschi}
\address{Giovanni A. Bonaschi \\ Institute for Complex Molecular Systems and Department of Mathematics and Computer Science, Technische Universiteit Eindhoven, P.O. Box 513, 5600 MB, Eindhoven, The Netherlands \&  Dipartimento di Matematica, Universit\`{a} di Pavia, 27100 Pavia, Italy} 
\email{g.a.bonaschi@tue.nl}

\author{J. A. Carrillo}
\address{Jos\'{e} A. Carrillo \\ Department of Mathematics, Imperial College London, London SW7 2AZ, United Kingdom}
\email{carrillo@imperial.ac.uk}

\author{M. Di Francesco}
\address{Marco Di Francesco\\Mathematical Sciences, University of Bath, Claverton Down, Bath BA2 7AY, United Kingdom}
\email{m.difrancesco@bath.ac.uk}

\author{M. A. Peletier}
\address{Mark A. Peletier\\ Institute for Complex Molecular Systems and Department of Mathematics and Computer Science, Technische Universiteit Eindhoven, Den Dolech 2, P.O. Box 513, 5600 MB Eindhoven, The Netherlands}
\email{M.A.Peletier@tue.nl}

\date{ \today}

\begin{document}

\begin{abstract}
We prove the equivalence between the notion of Wasserstein gradient flow for a one-dimensional nonlocal transport PDE with attractive/repulsive Newtonian potential on one side, and the notion of entropy solution of a Burgers-type scalar conservation law on the other. The solution of the former is obtained by spatially differentiating the solution of the latter. The proof uses an intermediate step, namely the $L^2$ gradient flow of the pseudo-inverse distribution function of the gradient flow solution. We use this equivalence to provide a rigorous particle-system approximation to the Wasserstein gradient flow, avoiding the regularization effect due to the singularity in the repulsive kernel. The abstract particle method relies on the so-called wave-front-tracking algorithm for scalar conservation laws. Finally, we provide a characterization of the sub-differential of the functional involved in the Wasserstein gradient flow.
\end{abstract}


\maketitle

\section{Introduction}
In this paper we construct and discuss connections between two partial differential equations on the real line.
The first equation is the nonlocal interaction equation
\begin{equation}\label{eq:main}
    \partial_t\mu  = \partial_x\bigl(\mu \;\partial_x W\mathord *\mu\bigr), \qquad x \in \R,\quad t>0,
\end{equation}
where $W$ is either the repulsive or the attractive Newton potential in one space dimension
\begin{equation}\label{eq:potentials}
    W(x)=-|x|\qquad \text{or} \qquad W(x)=|x|.
\end{equation}
We consider \emph{measure} solutions on the real line, with a given initial condition $\mu_0\in \PP$, where $\PP$ is the space of probability measures on $\R$ with finite second moment. The equation \eqref{eq:main} can be written as a continuity equation $\partial_t \mu + \partial_x (v\mu) = 0$ with $v:=-\partial_x W \ast \mu$. The velocity field $v(t,x)$ can be interpreted as the result of a \emph{nonlocal interaction} through the potential $W$, and the equation \eqref{eq:main} itself can be (at least formally) interpreted as the Wasserstein gradient flow of the following functional defined on $\PP$,
\begin{equation*}
  \W[\mu]=\frac{1}{2}\int_{\R \times \R}W(x-y)d\mu(x)d\mu(y).
\end{equation*}
The second equation is the scalar nonlinear conservation law
\begin{equation}\label{eq:main2}
    \partial_t F+ \partial_x g(F)=0, \qquad x \in \R,\;\; t>0,
\end{equation}
with
\begin{equation*}
  g(F)=F^2-F\qquad \hbox{or}\qquad g(F)=-F^2 +F.
\end{equation*}
We consider weak solutions $F$ on the real line, with a given initial condition $F_0$. The connection between the two equations \eqref{eq:main} and \eqref{eq:main2} is established through the relationship
\[
F(x) = F_\mu (x) := \mu\bigl((-\infty,x]\bigr).
\]
Formally, if $\mu$ is a solution of~\eqref{eq:main}, then $F$ is a solution of~\eqref{eq:main2}, and vice versa. This can be recognized by integrating~\eqref{eq:main} over $(-\infty,x]$ and observing that
\[
\partial_x (|\cdot| \mathord * \mu)=  (\sign\mathord *\mu) =  (1-2F).
\]
Because of this interpretation of $F$ in terms of $\mu$, we restrict ourselves to solutions $F$ of~\eqref{eq:main2}  that are increasing and bounded.

\medskip

These two equations have been studied extensively, but in different communities and using different tools. Equation \eqref{eq:main} describes the evolution of a system of interacting particles with an attractive ($W(x)=|x|$) or repulsive ($W(x)=-|x|$) potential, and equations of this type arise in a variety of physical, chemical, and biological applications; see e.g. \cite{golse,mogilner,TBL06} and the references therein. The specific example of the Newtonian potential in two dimensions arises in the Patlak-Keller-Segel model (see e.g.~\cite{Pat53,KS70,jager,BDP06}), where nonlocal transport effects are coupled with linear diffusion. Combined attractive/repulsive interactions have been studied in \cite{FR11,FR211,FHK11,CFP12,BCLR,BCLR2}. A deeper study on singular potentials has been performed in \cite{LT04,BD08,CDFLS,bertozzi1,bertozzi2,bertozzi3,bertozzi4}, see also the recent preprint \cite{CCH13}. In this paper, we focus on the two cases mentioned above, the attractive and the repulsive Newtonian potential in one dimension.

Similarly, equation~\eqref{eq:main2} has a long history and a wide range of applications; see e. g. \cite{bressan_book} and the references therein. Up to the unimportant linear term, the nonlinearity $g$ corresponds to the inviscid convex or concave Burgers equation, or Whitham's forward or backward equation~\cite{Whi74} (depending on the sign of $g$).

\medskip

For both equations, well-posedness strongly depends on the choice of the solution concept. Burgers' equation~\eqref{eq:main2} admits weak $L^\infty$ solutions on $\R$ for both choices of $g$~\cite[Sec.~3.4]{Eva98}, but there are examples of $L^\infty$ initial data which produce more than one solution. In order to single out physically relevant solutions in the context of gas dynamics, Oleinik \cite{Ole63} and Kru\v{z}kov \cite{Kru69,Kru70} formulated the concept of \emph{entropy solution}, which can be reached e.g.\ via a vanishing-viscosity approximation (see~\cite{bressan_book}). Different approximations give rise to other types of solutions, with so-called non-classical shocks~\cite{Lef02,VPP07}.

In the case of equation~\eqref{eq:main}, when $W$ is smooth and satisfies suitable growth bounds, distributional solutions exist and are unique. This follows as a trivial consequence of the theory in \cite{AGS}, but it can be easily deduced from minor modifications of the arguments in \cite{dobrushin}. For a less regular $W$, a distributional definition of a solution may not be meaningful; for instance, whenever $W'$ is discontinuous, the product $\delta (\partial_x W\mathord * \delta) = \delta W'$ is not well-defined. In this case the theory of Wasserstein gradient flows ~\cite{AGS,CDFLS} provides a solution concept for which existence and uniqueness holds provided $W$ is $\lambda$-convex, i.e. convex up to a quadratic perturbation. Therefore, the case $W(x)=|x|$ can be easily covered in view of the convexity of $W$ (see~\cite{CDFLS}); although $W(x)=-|x|$ is neither convex nor $\lambda$-convex, the corresponding functional (see~\eqref{eq:functional} below) \emph{is} $\lambda$-convex in the sense of McCann \cite{mccann}, see ~\cite{BS11,CFP12}, and therefore the abstract Wasserstein gradient flow theory applies.

\medskip

\begin{figure}[b]
\centering
\begin{tabular}{c | c | c | c |}
& Wasserstein gradient flow& Entropy solution & $L^2$ gradient flow\\
\hline & & & \\
\begin{tikzpicture}[scale=0.62]
\draw (0,3) node[rotate=90] {Attractive};
\draw (0,0) {};
\end{tikzpicture}
&
\begin{tikzpicture}[scale=0.7]
\draw (-3.3,0) node[anchor=south] {};
\draw (3.3,0) node[anchor=south] {};
\draw (0,-.5) node[anchor=south] {};
\draw [semithick, color=black] (0,-.2) -- (0,4.6);
\draw [semithick, color=black] (-3,0) -- (3,0);
\draw [semithick, color=black] (-2.5,0) -- (-2.5,1);
\draw [semithick, color=black] (-2.5,1) -- (-1.2,1);
\draw [semithick, color=black] (-1.2,1) -- (-1.2,0);
\draw [semithick, color=black] (2.5,0) -- (2.5,1);
\draw [semithick, color=black] (2.5,1) -- (1.2,1);
\draw [semithick, color=black] (1.2,1) -- (1.2,0);
\draw [dashdotted, color=black] (2,0) -- (2,1.5);
\draw [dashdotted, color=black] (2,1.5) -- (1,1.5);
\draw [dashdotted, color=black] (1,1.5) -- (1,0);
\draw [dashdotted, color=black] (-2,0) -- (-2,1.5);
\draw [dashdotted, color=black] (-2,1.5) -- (-1,1.5);
\draw [dashdotted, color=black] (-1,1.5) -- (-1,0);
\draw [->, dotted, color=black] (-0.5,0) -- (-0.5,2.3);
\draw [->, dotted, color=black] (0.5,0) -- (0.5,2.3);
\draw [->, ultra thick, color=black] (0,0) -- (0,4);
\draw[->, thin, color=black] (-2.3,.5) to [out=20,in=250] (-.3,2.1);
\draw[->, thin, color=black] (2.3,.5) to [out=160,in=290] (.3,2.1);
\end{tikzpicture}
&
\begin{tikzpicture}[scale=0.7]
\draw (-3.3,0) node[anchor=south] {};
\draw (3.3,0) node[anchor=south] {};
\draw (0,-3.5) node[anchor=south] {};
\draw [semithick, color=black] (0,-.2) -- (0,1.6);
\draw [semithick, color=black] (-3,0) -- (3,0);
\draw [->, semithick, color=black] (0,-0.4) -- (0,-1.2);
\draw [semithick, color=black] (0,-3.2) -- (0,-1.4);
\draw [semithick, color=black] (-3,-3) -- (3,-3);
\draw [thick, color=black] (-3,0) -- (-2,0);
\draw [thick, color=black] (2,1) -- (3,1);
\draw [thick, color=black] (-1,.5) -- (1,.5);
\draw [thick, color=black] (-2,0) -- (-1,.5);
\draw [thick, color=black] (1,.5) -- (2,1);
\draw [thick, color=black] (-.1,1) -- (.1,1);
\draw [dashdotted, color=black] (-1.3,0) -- (-.7,.5);
\draw [dashdotted, color=black] (1.3,1) -- (.7,.5);
\draw [dashdotted, color=black] (1.3,1) -- (2,1);
\draw [dotted, color=black] (-.4,0) -- (-.4,.5);
\draw [dotted, color=black] (.4,.5) -- (.4,1);
\draw [dotted, color=black] (.4,1) -- (1.3,1);
\draw [ultra thick, color=black] (-3,-3) -- (0,-3);
\draw [ultra thick, color=black]  (0,-2) -- (3,-2) ;
\fill[color=black] (0,-2) circle (0.5ex);
\draw[->, thin, color=black] (1.3,.75) -- (.6,.75);
\draw[->, thin, color=black] (-1.3,.25) -- (-.6,.25);
\end{tikzpicture}
&
\begin{tikzpicture}[scale=0.7]
\draw (-.5,0) node[anchor=south] {};
\draw (5.3,0) node[anchor=south] {};
\draw (0,-.5) node[anchor=south] {};
\draw [semithick, color=black] (0,-.2) -- (0,4.6);
\draw [semithick, color=black] (-.2,0) -- (5,0);
\draw [semithick, color=black] (4.5,-.1) -- (4.5,.1);
\draw [ultra thick, color=black] (0,2) -- (4.5,2);
\draw [thick, color=black] (0,.2) -- (2.25,1.1);
\draw [thick, color=black] (4.5,3.8) -- (2.25,2.9);
\draw [dashdotted, color=black] (0,.9) -- (2.25,1.4);
\draw [dashdotted, color=black] (4.5,3.1) -- (2.25,2.6);
\draw [dotted, color=black] (0,1.7) -- (2.25,1.7);
\draw [dotted, color=black] (2.25,2.3) -- (4.5,2.3);
\draw[->, thin, color=black] (1,.8) -- (1,1.6);
\draw[->, thin, color=black] (3.5,3.2) -- (3.5,2.4);
\end{tikzpicture}
\\
\hline & & & \\
\begin{tikzpicture}[scale=0.7]
\draw (0,2.7) node[rotate=90] {Repulsive};
\draw (0,0) {};
\end{tikzpicture}
&
\begin{tikzpicture}[scale=0.7]
\draw (-3.3,0) node[anchor=south] {};
\draw (3.3,0) node[anchor=south] {};
\draw (0,-.5) node[anchor=south] {};
\draw [->, ultra thick, color=black] (0,0) -- (0,4);
\draw [semithick, color=black] (0,-.2) -- (0,4.6);
\draw [semithick, color=black] (-3,0) -- (3,0);
\draw [dotted, color=black] (-.5,0) -- (-.5,3.4);
\draw [dotted, color=black] (-.5,3.4) -- (.5,3.4);
\draw [dotted, color=black] (.5,3.4) -- (.5,0);
\draw [dashdotted, color=black] (-1.2,0) -- (-1.2,2);
\draw [dashdotted, color=black] (-1.2,2) -- (1.2,2);
\draw [dashdotted, color=black] (1.2,2) -- (1.2,0);
\draw [thick, color=black] (-2.3,0) -- (-2.3,1);
\draw [thick, color=black] (-2.3,1) -- (2.3,1);
\draw [thick, color=black] (2.3,1) -- (2.3,0);
\draw[->, thin, color=black] (-.3,2.3) to [out=250,in=20] (-2.1,.5);
\draw[->, thin, color=black]  (.3,2.3) to [out=290,in=160] (2.1,.5);
\end{tikzpicture}
&
\begin{tikzpicture}[scale=0.7]
\draw (-3.3,0) node[anchor=south] {};
\draw (3.3,0) node[anchor=south] {};
\draw (0,-3.5) node[anchor=south] {};
\draw [semithick, color=black] (0,-.2) -- (0,1.6);
\draw [semithick, color=black] (-3,0) -- (3,0);
\draw [->, semithick, color=black] (0,-0.4) -- (0,-1.2);
\draw [semithick, color=black] (0,-3.2) -- (0,-1.4);
\draw [semithick, color=black] (-3,-3) -- (3,-3);
\draw [semithick, color=black] (-.1,-2) -- (.1,-2);
\draw [ultra thick, color=black] (-3,0) -- (0,0);
\draw [ultra thick, color=black]  (0,1) -- (3,1) ;
\draw [thick, color=black] (-2,-3) -- (2,-2);
\draw [thick, color=black] (-3,-3) -- (-2,-3);
\draw [thick, color=black] (2,-2) -- (3,-2);
\draw [dashdotted, color=black] (-1,-3) -- (1,-2);
\draw [dashdotted, color=black] (1,-2) -- (2,-2);
\draw [dotted, color=black] (-0.4,-3) -- (0.4,-2);
\draw [dotted, color=black] (0.4,-2) -- (1,-2);
\fill[color=black] (0,1) circle (0.5ex);
\draw[<-, thin, color=black] (1.5,-2.25) -- (.3,-2.25);
\draw[<-, thin, color=black] (-1.5,-2.75) -- (-.3,-2.75);
\end{tikzpicture}
&
\begin{tikzpicture}[scale=0.7]
\draw (-.5,0) node[anchor=south] {};
\draw (5.3,0) node[anchor=south] {};
\draw (0,-.5) node[anchor=south] {};
\draw [semithick, color=black] (0,-.2) -- (0,4.6);
\draw [semithick, color=black] (4.5,-.1) -- (4.5,.1);
\draw [semithick, color=black] (-.2,0) -- (5,0);
\draw [ultra thick, color=black] (0,2) -- (4.5,2);
\draw [dotted, color=black] (0,1.5) -- (4.5,2.5);
\draw [dashdotted, color=black] (0,1) -- (4.5,3);
\draw [thick, color=black] (0,.5) -- (4.5,3.5);
\draw[<-, thin, color=black] (1,.8) -- (1,1.6);
\draw[<-, thin, color=black] (3.5,3.2) -- (3.5,2.4);
\end{tikzpicture}
\\
\hline
\end{tabular}
\caption{Three forms of the same solution, for the gradient flow definition of equation~\eqref{eq:main} (left), for the entropy solution of~\eqref{eq:main2} (middle), and the $L^2$ gradient flow (right, Section~\ref{subsec:L2}). The top row is for the attractive case, the bottom row for the repulsive case. The direction of the evolution is indicated by arrows; the vertical arrows in the left column are Dirac delta functions. Note how the top and bottom evolutions are \emph{not} each other's time reversal}
\label{fig:schematicpic}
\end{figure}
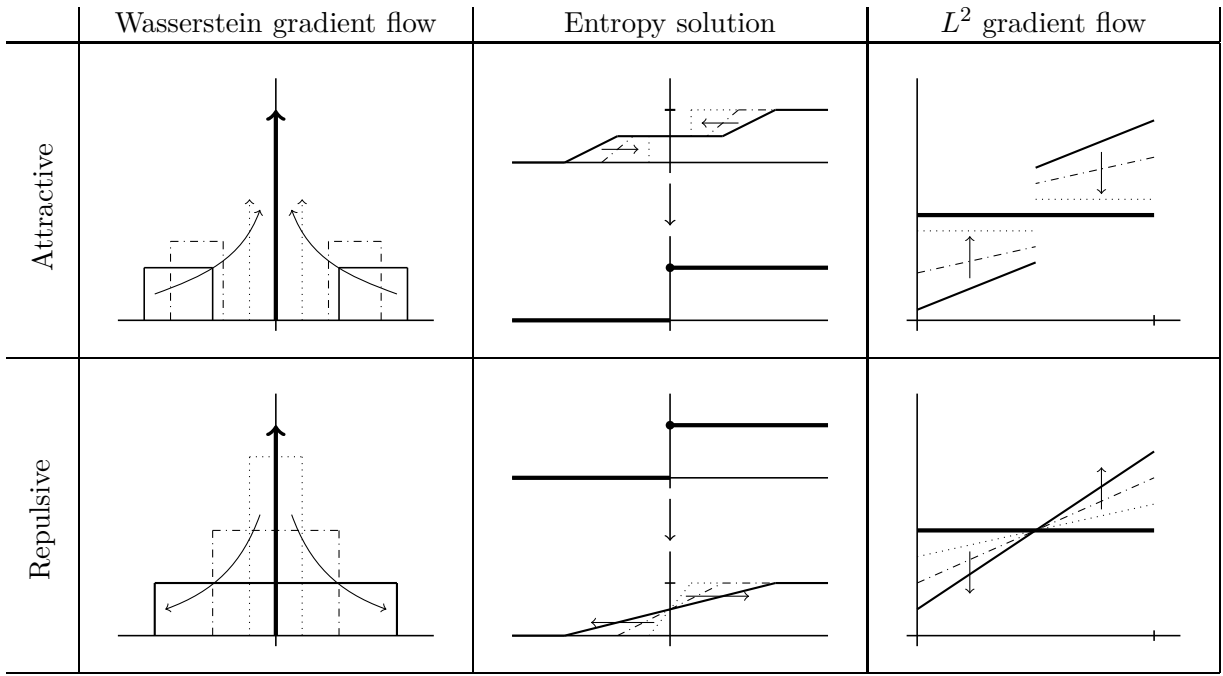

These well-posedness issues are strongly connected with the behaviour of the equations under time reversal. In Figure~\ref{fig:schematicpic} we illustrate this with an example. In the first column, solutions of~\eqref{eq:main} in the gradient flow concept are shown, with both attractive and repulsive interaction. In the attractive case, the two square waves collapse in finite time into Dirac delta functions and then propagate until they aggregate into a single delta function, which is a stationary solution. For the repulsive case, however, a single delta function is not stationary: it immediately regularizes into a square wave with linearly expanding boundaries. This example shows how the attractive and repulsive evolutions are not each other's time reversal.

The corresponding solutions of the Burgers equation are shown in the second column. The initial aggregation into delta functions translates into the formation of two shocks, which subsequently aggregate into a single fixed shock. With the opposite sign, the entropy condition disallows the corresponding time-reversed solution, and a rarefaction wave is formed instead.
These features are not limited to these examples; they occur for very general classes of initial data.

\medskip

Both the well-posedness subtleties and the non-invariance under time reversal raise questions about the connection between the two problems. For instance, how does the non-uniqueness in Burgers' equation manifest itself after transforming to~\eqref{eq:main}? What form does an entropy condition such as Oleinik's (see Definition~\ref{def:entropy}) take for solutions of~\eqref{eq:main}? Why does the gradient flow theory provide uniqueness for solutions of~\eqref{eq:main}, without further conditions? And is the unique gradient flow solution of~\eqref{eq:main} the same as the entropy solution of~\eqref{eq:main2}?

\subsection{Results}\label{subsec:results}

In the rest of this paper we address the above questions. Our main results are as follows.

\medskip

\begin{figure}[h]
\centering
\begin{tikzpicture}[scale=0.8]
\draw (0,0) node[anchor=center] {\Large $\mu_0$};
\draw [<->, thick, color=black] (.9,0) -- (3.1,0);
\draw (4,0) node[anchor=center] {\Large $F_0$};
\draw [<->, thick, color=black] (4.9,0) -- (7.1,0);
\draw (8,0) node[anchor=center] {\Large  $X_0$};
\draw [->, thick, color=black] (8,-.5) -- (8,-1.5);
\draw (4,-2) node[anchor=center] {\Large  $\widetilde{F}_t$};
\draw (4,-2.5) node[anchor=center] {\scriptsize $\parallel$};
\draw (3.4,-2.5) node[anchor=center] {\mycirc A};
\draw [->, thick, color=black] (4,-.5) -- (4,-1.5);
\draw (8,-2) node[anchor=center] {\Large  $\widetilde{X}_t$};
\draw (8,-2.5) node[anchor=center] {\scriptsize $\parallel$};
\draw (7.4,-2.5) node[anchor=center] {\mycirc B};
\draw [->, thick, color=black] (0,-.5) -- (0,-2.5);
\draw (0,-3) node[anchor=center] {\Large  $\mu_t$};
\draw [<->, thick, color=black] (.9,-3) -- (3.1,-3);
\draw (4,-3) node[anchor=center] {\Large  $F_t$};
\draw [<->, thick, color=black] (4.9,-3) -- (7.1,-3);
\draw (8,-3) node[anchor=center] {\Large  $X_t$};
\draw [<->, semithick, color=black] (-4,-1) -- (-2.5,-1);
\draw (-3.2,-.6) node[anchor=center] {\scriptsize define};
\draw [->, semithick, color=black] (-3.9,-1.5) -- (-3.9,-2.9);
\draw (-3.1,-2) node[anchor=center] {\scriptsize evolve};
\draw (-3.1,-2.4) node[anchor=center] {\scriptsize into};

\end{tikzpicture}
\caption{This figure illustrates the content of Theorem~\ref{thm:equivalence}. For a given initial datum, expressed in terms of $\mu_0$, $F_0$, or $X_0$, the theorem states that the three solutions at later time $t>0$ also are equivalent. Our contribution is the two equivalences A and~B}
\label{sketchproof}
\end{figure}
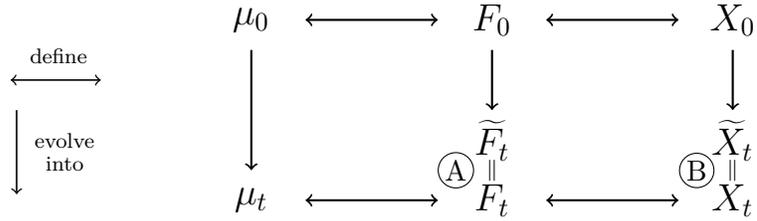

In Theorem \ref{thm:equivalence} we show that the gradient flow solution concept of~\eqref{eq:main} is equivalent to the entropy solution concept for~\eqref{eq:main2}. We establish this equivalence through a third solution concept, the $L^2$ gradient flow for the pseudo-inverse function $X$, which is defined in terms of $\mu$ and $F$ by
\[
X_{\mu}(s):= \inf \{ x \, | F_{\mu}(x) > s \}, \qquad s \in (0,1) ,
\]
and it maps $(0,1)$ to the support of $\mu$ (see Section~\ref{subsec:L2}). The content of the theorem is illustrated graphically in Figure~\ref{sketchproof}.

The proof is achieved by an explicit calculation for the case when $\mu_0$ is a  sum of delta functions (and $F_0$ and $X_0$ therefore both piecewise constant); the general case follows using the contractivity of the semigroup.

This result is the core of this paper: equations~\eqref{eq:main} and~\eqref{eq:main2} \emph{are} equivalent, provided one takes the `right' solution concept for both. In the latter we will discuss in detail how the specific aspects of the gradient flow concept and the entropy solution concept tie together

Let us mention here that the results in \cite{BBL05} already pointed out a link between scalar conservation laws with monotone data and the $L^2$ gradient flow. Moreover, it is worthwhile recalling that similar links between gradient flow solutions and entropy solutions have been lately explored in several contexts, see e.g. \cite{GO,DFM}.

\medskip

An important difference between the attractive and the repulsive case arises when one tries to approximate \emph{continuum} solutions to the Wasserstein gradient flow \eqref{eq:main} with a system of \emph{interacting particles}. Such a system typically reads as follows
\begin{equation*}
  \dot{x}_j(t)=-\frac{1}{N}\sum_{k=1}^N W'(x_j(t)-x_k(t)),\qquad j=1,\ldots,N,
\end{equation*}
and the approximation property is typically stated as
\begin{equation*}
  \frac{1}{N}\sum_{j=1}^N \delta_{x_j(t)} \rightharpoonup \mu(t)\quad \hbox{as}\ \ N\rightarrow +\infty,\qquad \hbox{for a. e.}\ \ t>0,
\end{equation*}
where the limit is intended in the weak-$\ast$ sense of measures, and $\mu(t)$ is the gradient flow solution to \eqref{eq:main}. When $W$ is smooth, say $C^2$, the above approximation property is easily recovered as delta type solutions $\frac{1}{N}\sum_{j=1}^N \delta_{x_j(t)}$ turn out to be a special case of gradient flow solutions; such a property is stated in short by saying that \emph{particles remain particles} in \eqref{eq:main}. As we already pointed out before, such a property may not be satisfied in case of a discontinuous $W$, since particle solutions may not be well defined because of the singularity in the \emph{self-interaction} force term $W'(0)$. Let us now focus on our case \eqref{eq:potentials}. In the attractive case $W(x)=|x|$ the results in \cite{CDFLS} provide a simple answer: \emph{particles remain particles, with the convention that the self-interaction term is neglected}. This is not surprising, as the force field is attractive, and e. g. two particles are not expected to exert forces on each other once they have collided. Let us mention that the result in \cite{CDFLS} holds in arbitrary dimension. In the repulsive case, the situation is way less trivial. This is already quite clear from the time reversal argument above: one single particle subject to the \emph{self-repulsive} force generates a squared (continuum) wave, and therefore it is clear that, in general, particles do not remain particles.

In Theorem \ref{thm:particle} we prove that a discrete approximation scheme for the Wasserstein gradient flow in the repulsive case can be constructed by exploiting the equivalence of \eqref{eq:main} with the scalar conservation law \eqref{eq:main2}. The appoximating procedure is based on the so-called \emph{wave-front-tracking method} (WFT) for conservation laws, see \cite{Daf72,dip76,Bre92}. This method consists mainly of two ingredients: discretization of initial data ($F_0 \to F_0^N$) and piecewise linear interpolation of the flux ($g \to g^N$). The peculiar characteristic is the discretization of the flux and that the two procedures are intimately related. This prevents the evolution from immediately regularizing any initial shock into a rarefaction wave. For the sake of completeness, we show that, for every positive time, the solution given by the WFT method is an approximation of the original solution. The proof is actually much simpler in our case, and it does not require the usual machinery used in the general theory for scalar conservation laws. Then, thanks to the equivalence result, we can rephrase such result into a particle approximation for the solution of the Wasserstein gradient flow. The final outcome is that, as in the attractive case, the self-repulsive force has to be neglected in the particle scheme. We point out that our result partially complements the results in the recent preprint \cite{CCH13}, in which a more general multi-dimensional theory is presented which does not cover the case of Newtonian potentials.

\medskip

During this work, a purely mathematical problem related to the definition of the Wasserstein sub-differential of $\W$ on singular measures came out, which is strictly related to the time reversal issue stated above. Collecting together the results from \cite{AGS,CDFLS,BS11}, one can prove existence and explicit characterization of the sub-differential of the functional $\W$ in the case of absolutely continuous measures for both the repulsive and the attractive case, and for concentrated measures in the attractive case. Unfortunately the same arguments cannot be applied when dealing with concentrated measures in the repulsive case. To handle this case we must refer to the more general (but less intuitive) notion of extended sub-differential (Definition \ref{def:frechet2}). Our analysis leads to the result in Proposition \ref{subdiffp}, which is an interesting example of extended sub-differential, with a geometrical view as well as with an explicit characterization. Two main properties are used in the proof: the $\lambda$-convexity of the functional and a closure property of the sub-differential.

\bigskip

The paper is organized as follows. In Section \ref{sec:concepts}, we introduce the three systems, with a particular attention at the Wasserstein and $L^2$ gradient flows where some results must be proven. For the part regarding entropy solutions we mainly refer to \cite{Eva98}. Section \ref{sec:equivalence} is devoted to rigorously prove the equivalence between the three concepts of solution. Section \ref{sec:particle} shows the applicability of the particle approximation. We finally study in Section \ref{sec:subdiff} the the sub-differential of $\W$ in detail, and give a characterization of its minimal element. We conclude with further discussion of the results of this paper in Section \ref{sec:discussion}.

\section{Three concepts of solutions}\label{sec:concepts}

In this section we give a precise definition of three solution concepts which we will show later on to be equivalent:
\begin{itemize}
  \item [(A)] Wasserstein gradient flow solution for \eqref{eq:main} (see Subsection \ref{subsec:wass})
  \item [(B)] Entropy solution for \eqref{eq:main2} (see Subsection \ref{subsec:entropy})
  \item [(C)] $L^2$ gradient flow for the pseudo-inverse equation obtained from \eqref{eq:main} (see Subsection \ref{subsec:L2})
\end{itemize}
For each of these notions we shall recall the existence and uniqueness results present in the literature, and complement them with some qualitative properties. The equivalence among the three notions will be proven rigorously in Section \ref{sec:equivalence}, and is supported here only by formal arguments. We stress here that the equivalence between (B) and (C) was suggested by the contractivity results in the Wasserstein distances for scalar conservation laws with monotone data proven in \cite{BBL05}. Our main contribution here is the link with the nonlocal interaction equation \eqref{eq:main} which was not described before.

\subsection{Wasserstein Gradient Flows}\label{subsec:wass}

Our starting point is that of the Wasserstein gradient flow in the space of probability measures in the spirit of \cite{AGS} combined with the recent results from~\cite{CDFLS,BS11,CFP12}. In what follows, $\PP$ is the space of probability measures on $\R$ with finite second moment. On the metric space $\PP$ endowed with the $2$-Wasserstein distance, we introduce the interaction energy functional
\begin{equation}\label{eq:functional}
  \W[\mu]=\frac{1}{2}\int_{\R \times \R}W(x-y)d\mu(x)d\mu(y),\qquad W(x)=\sigma|x|,\qquad \sigma\in \{-1,1\}.
\end{equation}
Next we recall the basic ingredients needed to define the notion of Wasserstein gradient flow, see~\cite{AGS}. First we define the push-forward measure. Let $\nu_1\in \mathcal{P}_2(\R^n)$ and let $T:\R^n\rightarrow \R^m$ be a $\nu_1$-measurable map. Then the push-forward measure of $\nu_1$ via $T$, denoted by $\nu_2=T_\sharp \nu_1\in \mathcal{P}_2(\R^m)$, is defined via $\nu_2(A)=\nu_1(T^{-1}A)$. For $i=1,2$ we recall the definition of $i$-th projection $\pi_i:\R\times \R\rightarrow \R$, $\pi_i(x_1,x_2)=x_i\in \R$. Given two measures $\mu_1,\mu_2\in \PP$, the $2$-Wasserstein distance between $\mu_1$ and $\mu_2$ is defined as following
\begin{equation*}
    d_W^2(\mu_1,\mu_2)=\min\left\{\int_{\R\times \R}|x-y|^2 d\bm{\gamma}(x,y)\;\; | \;\; \bm{\gamma} \in \mathcal{P}_2(\R\times \R),\;\; (\pi_i)_\sharp \bm{\gamma} =\mu_i,\: i=1,2\right\}.
\end{equation*}
The set of $\bm{\gamma} \in \mathcal{P}_2(\R\times \R)$ such that $(\pi_i)_\sharp \bm{\gamma} =\mu_i$ is called the set of \emph{plans} between $\mu_1$ and $\mu_2$, and is denoted by $\Gamma(\mu_1,\mu_2)$. The set of \emph{optimal plans} $\Gamma_0(\mu_1,\mu_2)\subset\Gamma(\mu_1,\mu_2)$ is the set of plans for which the minimum above is achieved, i.e. $\bm{\gamma}\in \Gamma_0(\mu_1,\mu_2)$ if and only if
\begin{equation*}
     d_W^2(\mu_1,\mu_2) = \int_{\R\times \R}|x-y|^2 d\bm{\gamma}(x,y).
\end{equation*}
Let $\mu_t \in \Abs([0,+\infty);\PP)$. The metric derivative of $\mu_t$ (if it exists) is given by
\begin{equation*}
    |\mu_t'|(t)=\lim_{h\rightarrow 0}\frac{d_W(\mu_{t+h},\mu_t)}{|h|}.
\end{equation*}
The metric derivative of an absolutely continuous curve is almost everywhere well defined, see~\cite{AGS}.

\begin{definition}[Fr\'{e}chet sub-differential]\label{def:frechet1}
Let $\phi:\PP\rightarrow (-\infty,+\infty]$ be proper and lower semi continuous, and let $\mu\in D(\phi)$. We say that $v \in L^2(\mu)$ belongs to the \emph{Frech\'{e}t sub-differential}, denoted by $\partial\phi(\mu)$, if
\begin{equation*}
    \phi(\mutil)-\phi(\mu)\geq \inf_{\bm{\gamma}\in \Gamma_0(\mu,\mutil)}\int_{\R\times \R} v(x)(y-x) d\bm{\gamma}(x,y)+ o(d_W(\mu,\mutil)).
\end{equation*}
\end{definition}

For $\mu\in \PP$ with $\partial \phi(\mu)\neq \emptyset$ we indicate $\partial^0 \phi(\mu)$ the element in $\partial \phi(\mu)$ with minimal $L^2(\mu)$-norm, which we refer to as the \emph{minimal sub-differential} of $\phi$ at $\mu$.
In some cases, this definition of sub-differential is too restrictive, and it should be replaced by the following one.

An important property needed to deal with Wasserstein gradient flow is $\lambda$-geodesic convexity of a functional. Let us first recall that, for $\mu,\nu\in \PP$, the curve $[0,1]\ni t \mapsto \mu_t=((1-t)\pi_1 + t\pi_2)_\sharp \bm{\gamma}$, with $\bm{\gamma}\in \Gamma_0(\mu,\nu)$, is a \emph{constant speed geodesic connecting $\mu$ to $\nu$}, i. e. it minimizes the action
\begin{equation*}
    \int_0^1 |\mu_t'|^2 dt,\qquad \hbox{on the set }\quad \mu_t\in \Abs([0,1];\PP)\quad\hbox{with}\;\;\mu_0=\mu,\;\; \mu_1=\nu.
\end{equation*}

\begin{definition}[$\lambda$-geodesic convexity]
Let $\phi:\PP\rightarrow (0-\infty,+\infty]$ be proper and lower semi-continuous, and let $\lambda\in \R$. Then, $\phi$ is $\lambda$-geodesically convex if, for all $\mu,\nu\in \PP$, there exists an optimal plan $\bm{\gamma}\in \Gamma_0(\mu,\nu)$ such that
\begin{equation*}
    \phi(\mu_t)\leq (1-t)\phi(\mu) + t\phi(\nu) -\frac{\lambda}{2}t(1-t)d_W^2(\mu,\nu),\qquad \hbox{for all}\;\; t\in [0,1],
\end{equation*}
where $\mu_t = ((1-t)\pi_1 + t\pi_2)_\sharp \bm{\gamma}$.
\end{definition}

Let us now turn back to our case, namely that of $\phi=\W$ in \eqref{eq:functional}. By combining the results in \cite{CDFLS,CFP12,BS11}, we obtain the following results. Here all the results are stated in one space dimension. We stress that the result in the following proposition in the attractive case is also valid in arbitrary space dimension.
\begin{proposition}\label{0geodesicconvex}
Let $W(x)$ be as in \eqref{eq:functional}. Then, the functional $\W$ is geodesically convex. Moreover, for all $\mu \in \PP$ such that it has no atoms (i.e. $\mu(\{ x\})=0$ for every $x \in \R$), the minimal Frech\'{e}t sub-differential $\partial^0\W(\mu)$ is well defined and contains the only element
\begin{equation}\label{eq:sub-differential_W_wass}
    \partial^0\W(\mu) =  \int_{x \neq y}\partial_x W(x-y)d\mu(y) = \sigma \int_{x \neq y}\sign(x-y) d\mu(y).
\end{equation}
Moreover, in the attractive case the formula \eqref{eq:sub-differential_W_wass} is valid for all $\mu\in \PP$.
\end{proposition}
The proof of the geodesically convexity relies on the representation of probability measures via pseudo-inverses of their distribution functions, and it will be proposed (in an equivalent form) in Proposition~\ref{convexityl2} in the next subsection. The characterization of the sub-differential in the general case of $\mu \in \PP$ is treated in Section \ref{sec:subdiff}.

\begin{definition}[Wasserstein Gradient flow]\label{def:grad_flow}
Let $W(x)$ as in \eqref{eq:functional} with $\sigma\in \{-1,1\}$. A curve
\[
\mu_t \in AC_{loc}^2([0,+\infty);\mathcal{P}_2(\R)),
\]
is a gradient flow for the functional $\W$ in \eqref{eq:functional} if there exists $v_t \in L^2(\mu_t)$ such that it satisfies
\begin{equation}\label{GradFlow}
\begin{split}
& \partial_t \mu_t + \partial_x(v_t\mu_t)=0 \quad \text{ in } \mathcal{D}'([0,+\infty) \times \R), \\
& v_t=-\partial^0\mathcal{W}[\mu_t] \quad  \text{ for a.e. }t>0.
\end{split}
\end{equation}
\end{definition}

The existence and uniqueness of gradient flow solutions in the sense of Definition \eqref{def:grad_flow} can be formulated in compact form, once again by combining the results in \cite{CDFLS,BS11}.

\begin{theorem}[Existence and uniqueness of gradient flows \cite{CDFLS,BS11}]\label{thm:dwcontraction}
Let $W(x)=\sigma|x|$ with $\sigma =\pm 1$ and $\mu_0\in \PP$. Then, there exists a unique (global-in-time) gradient flow solution for the functional $\W$ in the sense of Definition \eqref{def:grad_flow}, such that $\lim_{t \to 0}d_W(\mu_t,\mu_0)=0$. Moreover, for two given solutions $\nu_t$ and $\mu_t$, the following contraction property holds,
\begin{equation}\label{eq:contraction_wass}
    d_W(\nu_t,\mu_t)\leq d_W(\nu_0,\mu_0).
\end{equation}
Moreover, for $\sigma=-1$, the solution $\mu_t$ is absolutely continuous with respect to the Lebesgue measure for all $t>0$.
\end{theorem}

\begin{remark}
We could have stated the above definition by requiring $\partial^0\mathcal{W}[\mu_t]$ to be defined as in \eqref{eq:sub-differential_W_wass}. The main result in \cite{BS11} on the repulsive case implies in particular that $\mu_t$ is absolutely continuous respect to the Lebesgue measure, for all $t>0$ and for every inital $\mu_0$. Therefore, the explicit expression of the sub-differential can be used.
\end{remark}

\subsection{$L^2$ gradient flow}\label{subsec:L2}

Let us consider the Hilbert space $L^2((0,1))$ with norm $\| \cdot \|$, and the convex cone
\begin{equation*}
  \K:=\left\{f\in L^2((0,1))\, |\, f\ \hbox{ is non-decreasing}\right\}.
\end{equation*}
For a given $\mu\in \PP$, we define the cumulative distribution function $F_{\mu}(x)$ associated to $\mu$ as
\begin{equation*}
F_{\mu}(x):=\mu((-\infty,x]).
\end{equation*}
Then, we set $X_{\mu}$ as the pseudo-inverse of the the distribution function $F_{\mu}(x)$.
\begin{equation}
\label{pseudo-inverse}
X_{\mu}(s):=\inf \{ x: F_{\mu}(x)>s \} \qquad s \in (0,1).
\end{equation}
We can invert the above formula, and pass from $X_{\mu}$ to $F_{\mu}$, as follows
\begin{equation}\label{eq:rearrangement}
    F_{\mu}(x)=\int_0^1 \chi_{(-\infty,x]}(X_{\mu}(s))ds= \left| \left\{ X_{\mu}(s) \leq x \right\} \right|.
\end{equation}
In particular, both $F_\mu$ and $X_\mu$ are right-continuous and non-decreasing.
Now, given a probability measure $\mu \in \mathcal{P}(\R)$ and its pseudo-inverse $X_{\mu}$ we have
that
\begin{equation}\label{eq:changeofvar}
\int_{\mathbb{R}} \xi(x) d\mu(x)= \int_0^1\xi(X_{\mu}(s))ds,
\end{equation}
for every bounded continuous function $\xi$. Moreover, for $\mu,\nu\in \PP$, we can represent the Wasserstein distance $d_W(\mu,\nu)$ as
\begin{equation}\label{eq:distance_rappr}
d_W^2(\mu,\nu)=\int_0^1 \bigl|X_{\mu}(s)-X_{\nu}(s)\bigr|^2ds\,,
\end{equation}
and the optimal plan is given by $(X_{\mu}(s)\otimes X_{\nu}(s))_\# \mathcal{L}$, where $\mathcal{L}$ is the Lebesgue measure on the interval $[0,1]$.
These properties prove that there exists a natural isometry between $\PP$ and $\K\subset L^2([0,1])$, given by the mapping
\begin{equation*}
    \PP\ni \mu \mapsto X_\mu \in \K.
\end{equation*}
Through this identification it is possible to pose equation \eqref{eq:main} as a gradient flow in $L^2$ of a certain functional. In order to see that, let us first recall the following elementary computation already present in \cite{LT04,BS11}. Let $\mu_t$ be a gradient flow solution in the sense of Definition \ref{def:grad_flow} with no atoms for all times $t\geq 0$. Then, it is straightforward to find the following integro-differential equation satisfied by $X_t:=X_{\mu_t}$
\begin{equation}\label{eq:pseudo}
    \partial_t X_t(s) = -\sigma \int_0^1 \sign(X_t(s) -X_t(z)) dz,\qquad s\in [0,1],\quad t\geq 0.
\end{equation}
In order to give a meaning to \eqref{eq:pseudo} in case $X_t$ has atoms, we have to define $W'$ at zero. We assume henceforth that $W'(0)=0$.

In order to detect a gradient flow structure in $L^2$ for our equation \eqref{eq:main}, we should write $\W[\mu]$ in terms of the pseudo-inverse variable $X_{\mu_t}$. However, we have to make sure that the flow remains in the convex set $\K$. This procedure is reminiscent of \cite{Bre09}, see also \cite{BBL05}. Hence, the correct choice for the functional is the following. For a given $X\in L^2$, we set
\begin{align}
    & \mathscr{W}(X)= \frac{1}{2}\int_0^1\int_0^1 W(X(z)-X(\zeta))d\zeta dz, \nonumber\\
    & \mathscr{I}_\K(X) =
    \begin{cases}
    0 & \hbox{if}\quad X\in \K\\
    +\infty & \hbox{otherwise}
    \end{cases},    \label{eq:functional_l2} \\
    & \overline{\mathscr{W}}(X)= \mathscr{W}(X) + \mathscr{I}_\K(X). \nonumber \\
\end{align}

The functional $\mathscr{I}_\K$ is called the \emph{indicator function} of $\K$. Since the set $\K$ is convex, $\mathscr{I}_\K$ is a convex functional. We know that, for a given proper and lower semi-continuous functional $\mathscr{F}$ on $L^2((0,1))$, the sub-differential of $\mathscr{F}$ at $X\in L^2([0,1])$ is defined as the set
\begin{equation*}
    \partial \mathscr{F}(X)=\left\{ Y\in L^2((0,1))\;|\; \mathscr{F}(Z)-\mathscr{F}(X)\geq \int_0^1 Y(Z-X) + o(\|Z-X\|),\; \hbox{ for } \| Z-X \| \to 0 \right\}.
\end{equation*}

The sub-differential of the functional $\mathscr{I}_\K$ is characterized in the following proposition, which collects classical results in convex analysis plus more recent results from Brenier, Gangbo, Natile, Savar\'e and Westdickenberg \cite{NS09,BGSW13}. From now on, for a given element $X\in L^2({0,1})$, we use the notation
\begin{equation}\label{eq:omegaX}
  \Omega_X=\left\{s\in (0,1)\;|\;\; X\; \hbox{ is constant a.e. in a neighborhood of }\; s\right\},
\end{equation}
and note that $\Omega_X$ can always be written as a countable union of intervals, i.e. $\Omega_X = \bigcup_i I_i$.

\begin{proposition}[\cite{NS09,BGSW13}]\label{prop:natile}
Let $X\in \K$, and let $\Omega_X$ be defined as in \eqref{eq:omegaX}. Let
\begin{align*}
    & \mathcal{N}_X=\left\{ \mathcal{Z}\in C([0,1])\;\ |\;\;   \mathcal{Z}\geq 0\; \hbox{and}\;\;  \mathcal{Z}=0\;\;\hbox{in }\;\; [0,1]\setminus\Omega_X\right\}.
\end{align*}
For a given $Y\in L^2([0,1])$, let
\begin{equation*}
    \mathcal{Y}(s)=\int_0^s Y(\sigma)d\sigma.
\end{equation*}
Then, we have
\begin{equation*}
    Y\in \partial \mathscr{I}_\K(X)\qquad \Leftrightarrow\qquad  \mathcal{Y}\in \mathcal{N}_X.
\end{equation*}
In particular, if $Y\in \partial \mathscr{I}_\K(X)$, then
\begin{equation*}
    \begin{cases}
    Y=0 & \hbox{ a. e. in }\quad [0,1]\setminus \Omega_X\\
    \int_\alpha^\beta Y(s) ds =0 & \hbox{for every connected component $(\alpha,\beta)$ of $\Omega_X$}
    \end{cases}.
\end{equation*}
\end{proposition}

Let us now have a closer look at the functional $\mathscr{W}$. When restricted to $\K$, this functional can actually be proven to be \emph{linear}.
\begin{proposition}\label{prop:linear_functional}
Let $X\in \K$. Then
\begin{equation*}
    \mathscr{W}(X)= \sigma\int_0^1 (2z-1) X(z) dz.
\end{equation*}
\end{proposition}

\proof
We compute
\begin{align*}
    \mathscr{W}(X) &=\frac{\sigma}{2}\int_0^1\int_0^1|X(s)-X(z)|dz ds \\
    &= \frac{\sigma}{2}\int\int_{X(s)\geq X(z)}(X(s)-X(z))dz ds - \frac{\sigma}{2}\int\int_{X(s)\leq X(z)} (X(s)-X(z))dz ds\\
    &= \sigma\int\int_{X(s)\geq X(z)}(X(s)-X(z))dz ds,
\end{align*}
where we have used the symmetry of the two terms in the right hand side. Now, since $X$ is non-decreasing, the set $\{X(s)\geq X(z)\}$ can be written as
\begin{equation*}
    \{X(s)\geq X(z)\} = \{s\geq z\} \cup \{s\leq z \leq S(s)\},\qquad S(s)=\sup\{z\in [0,1]\;|\;\; X(z)=X(s)\},
\end{equation*}
and since $X(s)=X(z)$ on $\ \{s\leq z \leq S(s)\}$, we have
\begin{align*}
    \mathscr{W}(X) & =  \sigma\int\int_{s\geq z}(X(s)-X(z))dz ds = \sigma\left[ \int_0^1 \int_0^s X(s) dz ds -  \int_0^1 \int_z^1 X(z) ds dz\right] \\
    & = \sigma\left[ \int_0^1 s X(s) ds -  \int_0^1 (1-z) X(z) dz\right] = \sigma \int_0^1 X(z)(2z-1) dz.
\end{align*}
\endproof

An immediate consequence of Proposition \ref{prop:linear_functional} is the following
\begin{proposition}\label{convexityl2}
The functional $\overline{\mathscr{W}}$ is convex on $L^2([0,1])$.
\end{proposition}

\proof
Let $X_0, X_1\in L^2([0,1])$, and let $X_t=(1-t)X_0+tX_1$. If $X_0\in L^2([0,1])\setminus \K$, then $\overline{\mathscr{W}}(X_0)=+\infty$ and the inequality
\begin{equation*}
    \overline{\mathscr{W}}(X_t)\leq (1-t)\overline{\mathscr{W}}(X_0) + t\overline{\mathscr{W}}(X_1),
\end{equation*}
is trivially satisfied. The same holds if $X_1 \in L^2([0,1])\setminus \K$. On the other hand, if both $X_0,X_1\in \K$, then the above inequality is satisfied, since $\K$ is a convex set and $\mathscr{W}$ is linear.
\endproof

As another consequence of Proposition \ref{prop:linear_functional}, we have the following
\begin{proposition}\label{prop:gradient}
Let $X\in L^2([0,1])$. Then, $\partial \overline{\mathscr{W}}(X)\neq \emptyset$ if and only if $X\in \K$. In that case,
\begin{equation*}
    \partial \overline{\mathscr{W}}(X)\ni f(\cdot)\qquad f(s):=\sigma(2s-1),\quad s\in (0,1).
\end{equation*}
Moreover, if $X\in \K$ is strictly increasing, then $ \partial \overline{\mathscr{W}}(X)$ is single-valued and it therefore consists only of the $f$ defined above.
\end{proposition}

\proof
Assume $X\not\in \K$. Then, $\mathscr{I}_\K(X)=+\infty$. Hence, assuming the existence of $Y\in \partial \overline{\mathscr{W}}(X)$ implies
\begin{equation*}
    \mathscr{W}(Z)+\mathscr{I}_\K(Z) -\mathscr{W}(X) - \int_0^1 Y (Z-X)ds +o(\|X-Z\|) \geq \mathscr{I}_\K(X),
\end{equation*}
for all $Z\in L^2((0,1))$, i.e. in particular for all $Z\in \K$. But in the latter case, the left-hand side is finite whereas the right-hand side is infinite, which proves that $\partial \overline{\mathscr{W}}(X)=\emptyset$.

Let $X\in \K$ and $Z\in L^2((0,1))$. If $Z\not\in\K$, then the definition of sub-differential is trivially satisfied. Assume then $Z\in \K$ and by Proposition~\ref{prop:linear_functional}
\begin{equation*}
    \mathscr{W}(Z)-\mathscr{W}(X)=\sigma\int_0^1(2s-1)(Z(s)-X(s))ds.
\end{equation*}

Finally, assume that $X\in \K$ is strictly increasing. Suppose that there exists $g\in \partial \overline{\mathscr{W}}(X)$ with $g\neq f$ on an interval $I\subset (0,1)$. Let us assume without restriction that $g>f$ on $I$. Since $X$ is strictly increasing, there exists a $\bar X\in \K$ with $\bar X=X$ on $[0,1]\setminus I$ and $\bar X> X$ on $I$. Therefore, we have
\begin{align*}
    & \int_0^1 g(s)(\bar X - X) ds = \int_I  g(s)(\bar X - X) ds > \int_I  f(s)(\bar X - X) ds\\
    & \ = \int_0^1 f(s)(\bar X - X) ds = \mathscr{W}(\bar X)-\mathscr{W}(X),
\end{align*}
where the last step follows by Proposition \ref{prop:linear_functional}. Therefore, we have found an element $\bar X \in L^2((0,1))$ such that
\begin{equation*}
    \mathscr{W}(\bar X)-\mathscr{W}(X)< \int_0^1 g(s)(\bar X - X) ds,
\end{equation*}
and this contradicts the fact that $g\in \partial \overline{\mathscr{W}}(X)$. Therefore, $f$ is the only element in $\partial \overline{\mathscr{W}}(X)$.
\endproof

We now state the definition of gradient flow solution in $L^2$ for our problem.

\begin{definition}[$L^2$ gradient flow]\label{def_GFL2}
Let $W(x)=\sigma|x|$ with $\sigma\in\{-1,1\}$. An absolutely continuous curve $X_t \in L^2$ is an $L^2$ gradient flow for the functional $\overline{\mathscr{W}}$ defined in \eqref{eq:functional_l2} if it satisfies the differential inclusion
\begin{equation}\label{eq:GF_L2}
    - \partial_t X_t \in \partial \overline{\mathscr{W}}(X_t).
\end{equation}
\end{definition}

As $\overline{\mathscr{W}}$ is a convex functional on a Hilbert space, the classical theory of Brezis \cite{Bre73} can be applied to prove existence of a unique solution to \eqref{eq:GF_L2}.

\begin{theorem}[Existence and uniqueness of $L^2$ gradient flow]\label{thm:brezis}
Let $W(x)=\sigma|x|$ with $\sigma\in\{-1,1\}$ and let $X_0\in \K$. Then, there exists a unique gradient flow solution $X_t$ in the sense of Definition~\ref{def_GFL2} with initial condition $X_0$. Moreover, for two solutions $X_{0,t}$ and $X_{1,t}$ to \eqref{eq:GF_L2}, the following contraction property holds
\begin{equation}\label{eq:contraction_L2}
    \|X_{0,t} - X_{1,t}\| \leq  \|X_{0,0} - X_{1,0}\|,
\end{equation}
for all $t\geq 0$. Moreover, for $\sigma=-1$, the solution $X_t$ is strictly increasing.
\end{theorem}

As a byproduct of the theory in \cite{Bre73}, the \emph{minimal selection} of the sub-differential $\partial \overline{\mathscr{W}}$ is achieved in the differential inclusion \eqref{eq:GF_L2} at a.e. time. Since $\partial \overline{\mathscr{W}}(X)$ is a convex set, it admits a unique element of minimal norm, that we call the \emph{minimal sub-differential} of $\overline{\mathscr{W}}$ at $X$, and we denote by $\partial^0 \overline{\mathscr{W}}(X)$. We characterize the minimal sub-differential in both the attractive and the repulsive case in the following theorem. As the sub-differential is single valued in case $X$ is strictly increasing, clearly we shall restrict to the case $X\in \K$ such that $\Omega_X \neq \emptyset$. It must be noticed that the mathematical structure coincide perfectly with the Wasserstein framework, a strong sign revealing the equivalence.

\begin{theorem}\label{thm:minimal}
Let $X \in \K$ and let $\overline{\mathscr{W}}$ be as in \eqref{eq:functional_l2}. Let $\Omega_X=\bigcup_{j\in J} I_j$ with $J$ possibly empty, where $I_j=(\alpha_j,\beta_j)$ are ordered disjoint intervals. If $\sigma=-1$, then
\begin{equation*}
    \partial^0 \overline{\mathscr{W}}(X)(s) = -2s+1,\qquad \hbox{for all}\ \ s\in [0,1].
\end{equation*}
If $\sigma=1$, then
\begin{equation}\label{eq:minimal_L2_attr}
     \partial^0 \overline{\mathscr{W}}(X)(s) =
     \begin{cases}
     2s-1 & \hbox{if}\quad s\in [0,1]\setminus \Omega_X \\
     \alpha_j+\beta_j-1 & \hbox{if}\quad s\in I_j
     \end{cases}.
\end{equation}
\end{theorem}
The reader may be surprised of the term $\alpha_j+\beta_j-1$ that appears in the sub-differential when $\sigma=1$. In fact that term can be seen as $\frac{1}{2}(2\alpha_j-1)+\frac{1}{2}(2\beta_j-1)$, i.e. the average of the sub-differential evaluated at the two extrema of the interval.
\proof

Let $\sigma=-1$. By additivity of the sub-differential, all the elements $Y\in  \partial\overline{\mathscr{W}}(X)$ are of the form
\begin{equation*}
    Y(s)= -2s+1 + Z(s),
\end{equation*}
with $Z\in L^2([0,1])$ such that
\begin{equation*}
    \mathcal{Z}(s)=\int_0^s Z(\sigma)d\sigma,
\end{equation*}
satisfies $\mathcal{Z}\geq 0$ and $\mathcal{Z} =0$ in $[0,1]\setminus \Omega_X$. Now, let us compute
\begin{align*}
    & \|Y\|_{L^2}^2 = \int_0^1 (-2s+1 + Z(s))^2 ds = \int_0^1 (-2s+1)^2 ds + \int_0^1 Z(s)^2 ds + 2\int_0^1 (-2s+1)Z(s) ds\\
    & \ \ = \int_0^1 (-2s+1)^2 ds + \int_0^1 Z(s)^2 ds + \left[(-2s+1)\mathcal{Z}(s)\right]_{s=0}^{s=1} + 4\int_0^1\mathcal{Z}(s) ds.
\end{align*}
Now, since $s=0$ and $s=1$ are not elements in $\Omega_X$, clearly we have $\mathcal{Z}(0)=\mathcal{Z}(1)=0$. Therefore, the boundary term above vanishes. All the other terms are non-negative, and therefore the minimum of $\|Y\|_{L^2}^2$ is achieved with $Z\equiv 0$.

Assume now $\sigma=1$. Let us first check that $\partial^0 \overline{\mathscr{W}}(X)$ defined in \eqref{eq:minimal_L2_attr} belongs to $\partial\overline{\mathscr{W}}(X)$. We have to check that
\begin{equation*}
    \overline{\mathscr{W}}(\widetilde{X})-\overline{\mathscr{W}}(X)\geq \int_0^1 \partial^0 \overline{\mathscr{W}}(X)(s)(\widetilde{X}(s)-X(s)) ds.
\end{equation*}
Since the above inequality is trivially satisfied if $\widetilde{X}\not\in\K$, we can assume $\widetilde{X}\in \K$ and use Proposition \ref{prop:linear_functional}. We first assume $\widetilde{X}\in C^1$. We have to check
\begin{equation*}
    \int_0^1(2s-1)(\widetilde{X}(s)-X(s))ds\geq \int_0^1 \partial^0 \overline{\mathscr{W}}(X)(s)(\widetilde{X}(s)-X(s)) ds,
\end{equation*}
which, in view of \eqref{eq:minimal_L2_attr}, is equivalent to
\begin{equation}\label{eq:minimal_toprove}
    \sum_{j=1}^{\infty}\int_{\alpha_j}^{\beta_j} g_j(s)(\widetilde{X}(s)-x_j) ds\geq 0,
\end{equation}
where $x_i\equiv X|_{I_j}$, and defining $y_j=\alpha_j+\beta_j-1$
\begin{equation*}
    g_j(s):=2s-1 - y_j,\qquad \hbox{for}\ \ s\in (\alpha_j,\beta_j).
\end{equation*}
In order to prove \eqref{eq:minimal_toprove}, we first observe that
\begin{equation*}
    G_j(s):=\int_{\alpha_j}^s g_j(\sigma)d\sigma,\qquad \hbox{for}\ \ s\in (\alpha_j,\beta_j),
\end{equation*}
satisfies $G_j(s)\leq 0$ on $I_j$ and $G_j(\alpha_j)=G_j(\beta_j)=0$. Hence, since $\widetilde{X}\in C^1$, we can integrate by parts to obtain
\begin{align*}
    & \int_{\alpha_j}^{\beta_j} g_j(s)(\widetilde{X}(s)-x_j) ds = \left[ G_j(s) (\widetilde{X}(s)-x_j)\right]_{s=\alpha_j}^{s=\beta_j} - \int_{I_j}G_j(s)\widetilde{X}'(s) ds \\
    & \qquad = - \int_{I_j}G_j(s)\widetilde{X}'(s) ds \geq 0,
\end{align*}
since $\widetilde{X}'\geq 0$. The general case $\widetilde{X}\in \K$ can be easily obtained by approximation.

Now we have to check the minimality condition. As in the case $\sigma=-1$, we know that all the elements $Y\in \partial\overline{\mathscr{W}}(X)$ are of the form
\begin{equation*}
    Y(s)= 2s-1 + Z(s),
\end{equation*}
with the same conditions on $Z$ as in case $\sigma=-1$. Then, using the property of $Z$ in proposition \ref{prop:natile}, we get
\begin{align*}
    & \|Y\|_{L^2}^2 = \int_{[0,1]\setminus \Omega_X} (2s-1+Z(s))^2ds +\int_{\Omega_X} (2s-1+Z(s))^2ds \\
    & \qquad = \int_{[0,1]\setminus \Omega_X} (2s-1)^2ds +\int_{\Omega_X} (2s-1+Z(s))^2ds.
\end{align*}
Therefore, in order to achieve the minimal selection, we have to minimize
\begin{equation*}
    \int_{\Omega_X} (2s-1+Z(s))^2ds,
\end{equation*}
on the set of $Z\in L^2([0,1])$ such that
\begin{equation*}
    \mathcal{Z}(s)=\int_0^s Z(\sigma)d\sigma,
\end{equation*}
satisfies $\mathcal{Z}\geq 0$ and $\mathcal{Z} =0$ in $[0,1]\setminus \Omega_X$. Notice in particular that $Z$ has to satisfy the constraint $\int_{I_j}Z(s)ds =0$. Therefore, the minimal selection for $\int_{I_j}(2s-1+Z(s))^2ds$ should be sought in the class $\int_{I_j}(2s-1+Z(s))ds = y_j(\beta_j-\alpha_j)$. The previous equality holds because of the following formula: $y_j=\frac{1}{\beta_j-\alpha_j}\int_{I_j}(2s-1)ds$. A direct  argument in the minimization of the $L^2$ norm gives that the minimizer should be constant on $I_j$, with the constant being given by $y_j$. This gives
\begin{equation*}
    Z(s)=y_j -2s+1\qquad \hbox{on}\ \ I_j,
\end{equation*}
and the assertion is proven.
\endproof

The result in Theorem \ref{thm:minimal} allows to provide an explicit formula for the unique gradient flow solution provided in Theorem \ref{thm:brezis} in the repulsive case, and a more refined formula for the time derivative $\partial_t X_t$ in the attractive case. The proof is an elementary consequence of Theorem \ref{thm:minimal}, and is therefore omitted.

\begin{theorem}\label{thm:explicit}
Let $X_0\in \K$. If $\sigma=-1$, then, the unique gradient flow solution $X_t$ in the sense of Definition \ref{def_GFL2} with initial condition $X_0$ satisfies
\begin{equation}\label{eq:explicit1}
    X_t(s)=X_0(s) +t(2s-1),
\end{equation}
for all $s\in [0,1]$ and $t\geq 0$. If $\sigma=1$, given
\begin{equation*}
    \Omega_{X_t}=\bigcup_{j=1}^{+\infty}(\alpha_j(t),\beta_j(t)),
\end{equation*}
then $X_t$ satisfies
\begin{equation}\label{eq:explcit2}
     -\partial_t X_t(s) =
     \begin{cases}
     2s-1 & \hbox{if}\quad s\in [0,1]\setminus \Omega_{X_t} \\
     \alpha_j(t)+\beta_j(t)-1 & \hbox{if}\quad s\in (\alpha_j(t),\beta_j(t))
     \end{cases}.
\end{equation}
\end{theorem}

\subsection{Entropy solutions}\label{subsec:entropy}

We now turn our attention to the cumulative distribution variable
\begin{equation*}
F(x,t):=\mu_t((-\infty,x])=\int_{-\infty}^xd\mu_t(y),
\end{equation*}
where $\mu_t$ is a Wasserstein gradient flow in the sense of Definition \ref{def:grad_flow}. Assume for simplicity that $\mu_t=\rho\leb$, and that $\rho(\cdot,t)$ is compactly supported. Then,
\begin{align*}
    & \partial_t F = \int_{-\infty}^x \partial_t \rho_t (x) dx = \sigma F_x \int_{-\infty}^{+\infty} \sign(x-y) \rho (y,t) dy \\
    & \ \ = \sigma F_x\left(\int_{-\infty}^{x}\rho (y,t) dy - \int_{x}^{+\infty}\rho (y,t) dy\right) = \sigma F_x\left(2F(x,t)- 1\right) = \sigma \partial_x( F^2 - F),
\end{align*}
hence $F$ satisfies the \emph{scalar conservation law}
\begin{equation}\label{eq:CL}
   \partial_t F + \partial_x g(F) =0\quad \hbox{with }\ g(F)=\sigma F(1-F) \text{ and } \sigma \in \{ -1, 1 \}.
\end{equation}
As shock waves (discontinuities) may appear in finite time, a concept of weak solution is needed. As more than one weak solution may arise with the same initial condition, the concept of entropy solution \cite{Ole63} is needed, in order to select admissible shock waves.

\begin{definition}[Entropy solution]\label{def:entropy}
Let $g$ be  as in \eqref{eq:CL}, and let $F_0\in L^\infty(\R)$ be a non-decreasing function. A function $F \in L^{\infty}([0,+\infty)\times\R)$ is called \emph{entropy solution} if it is a solution of the following initial value problem
\begin{equation}\label{eq:entropic}
\begin{cases}
\partial_tF+ \partial_x g(F)=0 & \text{in }\mathcal{D}'((0,+\infty) \times \R) \\
F=F_0 & \text{on } \R \times (t=0)
\end{cases},
\end{equation}
and if, in the case $\sigma=1$ (i.e. $g$ convex), it satisfies the Oleinik condition:
\begin{equation}\label{Oleinik}
F(x+z,t)-F(x,t) \leq \frac{C}{t}z,
\end{equation}
for some constant $C\geq 0$ and a.e. $x \in \R$, $z,t > 0$.
\end{definition}

Notice that no Oleinik condition \cite{Ole63} is needed if $\sigma =1$, as decreasing (non entropic) jumps are excluded a-priori since our solutions are non decreasing. The existence and uniqueness of an entropy solution to \eqref{eq:entropic} is guaranteed by the classical result in \cite{Ole63}, see also \cite{Kru70}.

\begin{theorem}[Existence and uniqueness of $L^\infty$ entropy solutions]\label{thm:entropy}
Let $g$ be as in \eqref{eq:CL}, and let $F_0\in L^\infty(\R)$ non-decreasing. Then, there exists a unique entropy solution in the sense of Definition \ref{def:entropy} with initial condition $F_0$. Moreover, let $F_0,F_1 \in L^\infty(\R)$ with $F_0-F_1 \in L^1(\R)$. Then, the two entropy solutions $F_0(\cdot,t)$ and $F_1(\cdot,t)$ with initial conditions $F_0$ and $F_1$ respectively satisfy the contraction property
\begin{equation}\label{eq:contraction_entropy}
    \|F_1(\cdot,t)-F_2(\cdot,t)\|_{L^1(\R)}\leq  \|F_1(\cdot,0)-F_2(\cdot,0)\|_{L^1(\R)}.
\end{equation}
\end{theorem}

The contraction result of Theorem \ref{thm:entropy} was originally proven in \cite{Kru69}, and is well explained also in \cite[Proposition 2.3.6]{Ser99}.

\begin{remark}
Clearly, when $F_0$ is the cumulative distribution of a probability measure $F_0(x)=F_{\mu_0}(x)=\int_{-\infty}^x d\mu_0(x)$, then $F_0$ is non-decreasing on $\R$. It can be proven by means of classical results on the Burgers equation that $F(\cdot,t)$ is non decreasing for all times $t\geq 0$. More precisely, one can express the unique entropy solution via the Lax-Oleinik formula, cf. e.g. \cite[Section 3.4.2]{Eva98}, and use the monotonicity of $g'$ to prove the assertion. Since we will obtain the same property as a by-product of our results, we skip the details at this stage.
\end{remark}

For future use, we recall the notion of \emph{Riemann problem} for \eqref{eq:CL}. A Riemann problem is an initial value problem \eqref{eq:entropic} with initial condition
\begin{equation}\label{eq:riemann}
  F_0(x)=
  \begin{cases}
  F^L & \hbox{ if }\quad x<0 \\
   F^R & \hbox{ if }\quad x>0
  \end{cases},
\end{equation}
with $F^L< F^R$. The solution to the Riemann problem in this case depends on the sign of $\sigma$. If $\sigma=1$, then the flux $g$ is concave, therefore increasing shocks are admissible. On the other hand, if $\sigma=-1$, then the flux $g$ is convex, and increasing shock are not admissible, and the initial discontinuity in the Riemann problem is solved by a rarefaction wave. More precisely, the solution to \eqref{eq:riemann} in the case $\sigma=1$ is given by
\begin{equation*}
  F(x,t)=
  \begin{cases}
  F^L & \hbox{ if }\quad x<(1-(F^R+F^L))t \\
   F^R & \hbox{ if }\quad x>(1-(F^R+F^L))t
  \end{cases}.
\end{equation*}
We recall that the speed of propagation of the shock wave between $F^L$ and $F^R$ is obtained via the Rankine-Hugoniot condition
\begin{equation}\label{eq:RH}
    \dot{x}(t)=\frac{g(F^L)-g(F^R)}{F^L-F^R}.
\end{equation}
In the case $\sigma=-1$, the solution is given by
\begin{equation*}
  F(x,t)=
  \begin{cases}
  F^L & \hbox{ if }\quad x<(-1+2F^L)t \\
  \frac{x+t}{2t}& \hbox{ if }\quad (-1+2F^L)t<x<(-1+2F^R)t\\
   F^R & \hbox{ if }\quad x>(-1+2F^R)t
  \end{cases}.
\end{equation*}

\section{Equivalence of the three notions of solutions}\label{sec:equivalence}

The following theorem is the main result of this paper.

\begin{theorem}[Equivalence of the three solutions.]\label{thm:equivalence}
Let $W(x)=\sigma|x|$ with $\sigma \in\{-1,1\}$. Let $\mu_0 \in \mathcal{P}_2(\R)$. Let $F_0(x)=\mu_0((-\infty,x])$ and let $X_0$ be the pseudo-inverse of $F_0$. Let $g$ be defined as in \eqref{eq:CL}. Let $\mu_t\in \Abs([0,+\infty))\rightarrow \PP$ be any curve. Then, the following are equivalent:
 \begin{itemize}
   \item [(C1)] The curve $\mu_t$ is the unique gradient flow solution in the sense of Definition~\ref{def:grad_flow} with initial condition $\mu_0$.
   \item [(C2)] The curve $F(\cdot,t)=\mu_t((-\infty,x])$ is the unique entropy solution in the sense of Definition~\ref{def:entropy} with initial condition $F_0$.
   \item [(C3)] The curve $X_t(s)=\inf\{x| F(x,t)>s\}$ is the unique $L^2$ gradient flow in the sense of Definition~\ref{def_GFL2} with initial condition $X_0$.
 \end{itemize}
\end{theorem}

\proof

\textsc{Step 1 - Finite combination of delta measures}.

The proof is divided in two parts. In the first one we prove the equivalence only for initial conditions involving finite sum of delta measures, considering the attractive and the repulsive case separately. Then we prove the equivalence for any initial condition with an approximation argument.

\begin{figure}[h]
\centering
\begin{tikzpicture}
\draw[semithick, color=black] (-2,0) -- (-2,4);
\draw[semithick, color=black] (-2,0) -- (2,0);
\draw[semithick, color=black] (-2,-.1) -- (-2,0);
\draw (-2,-.1) node[anchor=north] {$0$};
\draw[semithick, color=black] (-1.5,-.1) -- (-1.5,.1);
\draw[semithick, color=black] (-1,-.1) -- (-1,.1);
\draw[semithick, color=black] (-.5,-.1) -- (-.5,.1);
\draw[semithick, color=black] (0,-.1) -- (0,.1);
\draw[semithick, color=black] (.5,-.1) -- (.5,.1);
\draw[semithick, color=black] (1,-.1) -- (1,.1);
\draw (1,-.1) node[anchor=north] {$1$};
\draw (1.9,0) node[anchor=north] {$s$};
\draw (-2,3.8) node[anchor=east] {$X_0(s)$};
\draw[semithick, color=black] (-2,1) -- (-1,1) circle (2pt);
\draw[semithick, color=black] (-2.1,1) -- (-1.9,1);
\fill[color=black] (-1,1.5) circle (0.3ex);
\draw[semithick, color=black] (-1,1.5)  -- (-.5,1.5) circle (2pt) ;
\draw[semithick, color=black] (-2.1,1.5) -- (-1.9,1.5);
\draw (-2,1.5) node[anchor=east] {$x_i$};
\fill[color=black] (-.5,2) circle (0.3ex);
\draw[semithick, color=black] (-.5,2) -- (.5,2) circle (2pt);
\draw[semithick, color=black] (-2.1,2) -- (-1.9,2);
\fill[color=black] (.5,3) circle (0.3ex);
\draw[semithick, color=black] (.5,3)  -- (1,3);
\draw[semithick, color=black] (-2.1,3) -- (-1.9,3);
\draw[semithick, color=black] (7,0) -- (7,4);
\draw[semithick, color=black] (3,0) -- (11,0);
\draw[semithick, color=black] (7,-.1) -- (7,0);
\draw (10.9,0) node[anchor=north] {$x$};
\draw (7,3.8) node[anchor=east] {$F_0(x)$};
\draw (7,3) node[anchor=east] {$1$};
\draw[semithick, color=black] (7.1,3) -- (6.9,3);
\draw[semithick, color=black] (7.1,.5) -- (6.9,.5);
\draw[semithick, color=black] (7.1,1) -- (6.9,1);
\draw[semithick, color=black] (7.1,2) -- (6.9,2);
\draw[semithick, color=black] (7.1,2.5) -- (6.9,2.5);
\draw[semithick, color=black] (5,1)  -- (6,1) circle (2pt);
\fill[color=black] (5,1) circle (0.3ex);
\draw[semithick, color=black] (5,.1) -- (5,-.1);
\draw[semithick, color=black] (6,1.5)  -- (7.5,1.5) circle (2pt);
\fill[color=black] (6,1.5) circle (0.3ex);
\draw[semithick, color=black] (6,.1) -- (6,-.1);
\draw (6,0) node[anchor=north] {$x_i$};
\draw[semithick, color=black] (7.5,2.5)  -- (8.5,2.5) circle (2pt);
\fill[color=black] (7.5,2.5) circle (0.3ex);
\draw[semithick, color=black] (7.5,.1) -- (7.5,-.1);
\draw[semithick, color=black] (8.5,3)  -- (11,3);
\fill[color=black] (8.5,3) circle (0.3ex);
\draw[semithick, color=black] (8.5,.1) -- (8.5,-.1);
\end{tikzpicture}
\caption{$X_0$ and $F_0$ corresponding to a concentrated $\mu_0$}
\label{fig:mu0f0}
\end{figure}
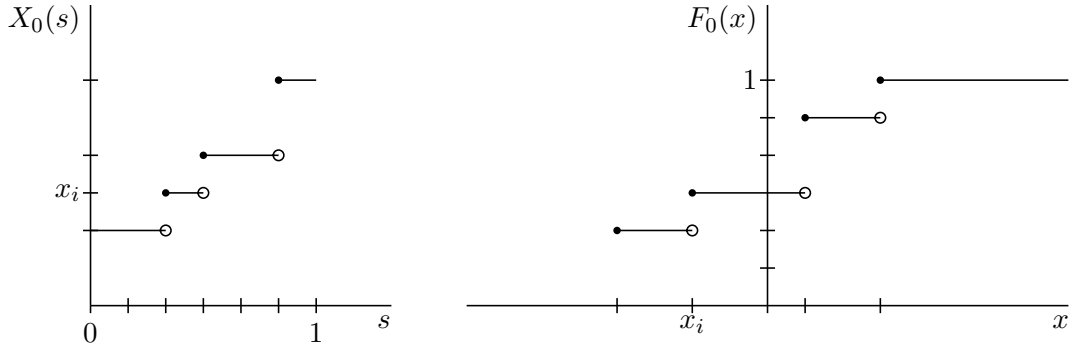

We first consider the class of initial conditions
\begin{equation}\label{eq:thm_deltas}
    \mu_0 = \sum_{j=1}^N m_j\delta_{x_j}, \qquad  1=\sum_{j=1}^N m_j.
\end{equation}
Let us set $M_0=0$ and $M_j=\sum_{k=1}^j m_k$, for $j=1,\ldots,N$. In particular, we have $M_N=1$. We easily get (see the example in Figure~\ref{fig:mu0f0})
\begin{equation*}
    F_0=\sum_{j=1}^{N} m_j \chi_{[x_j,+\infty)},\qquad X_0=\sum_{j=1}^{N} x_j \chi_{(M_{j-1},M_{j})}.
\end{equation*}

We now distinguish between the attractive and the repulsive case.

\emph{\textit{Attractive case}}. In the case $\sigma=1$, we claim that the unique $L^2$ gradient flow solution in the sense of definition \ref{def_GFL2} with initial condition $X_0$ is given by
\begin{equation*}
    X(s,t)=x_1(t)\chi_{[0,m_1)}+\sum_{j=2}^{N-1} x_j(t) \chi_{[M_{j-1},M_{j})}+ x_{N}(t)\chi_{[M_{N-1},1]},
\end{equation*}
with the $x_j$'s solving the particle system
\begin{align*}
    \dot{x}_j(t)=-\sum_{k:\ x_j(t)\neq x_k(t)} m_j\sign(x_j(t) - x_k(t)),
\end{align*}
with the convention that $\sign (0)=0$, so particles can collide and stick together. The proof of the claim is contained in Theorem \ref{thm:explicit}, checking that the velocity of the particles is $\dot{x}_j=1-\alpha_j-\beta_j$. Now, let $F(x,t)=\int_0^1 \mathcal{X}_{(-\infty,x]}(X(s,t))ds$. In \cite[Remark~2.10]{CDFLS} it is proven that the curve of probability measures $\mu_t=\partial_x F(x,t)$ is the unique Wasserstein gradient flow solution with initial condition $\mu_0$ in the sense of Definition \ref{def:grad_flow}. It remains to prove that $F(\cdot,t)$ is an entropy solution with initial condition $F_0$ in the sense of Definition \ref{def:entropy}. Let us first observe that in this case
\begin{equation*}
    F(x,t)=\mu_t((-\infty,x]) = \sum_{j=1}^{N} m_j \chi_{[x_j(t),+\infty)}.
\end{equation*}
Hence, we only need to prove that all the shocks in $F$ are admissible and that they satisfy the Rankine-Hugoniot condition \eqref{eq:RH}. To see this, let us compute
\begin{align*}
     & \dot{x}_j(t)= 1-M_j^-(t) -M_{j}^+(t),\\
     & M_j^-(t) :=\sum_{x_k(t)<x_j(t)} m_k,\qquad  M_j^+(t) :=\sum_{x_k(t) \leq x_j(t)} m_k.
\end{align*}
Clearly, the above identity yields
\begin{equation*}
    \dot{x}_j(t)= 1 - F(x_j(t)^-,t) - F(x_j(t)^+,t) = \frac{g(F(x_j(t)^-,t)) - g(F(x_j(t)^+,t))}{F(x_j(t)^-,t) - F(x_j(t)^+,t)},
\end{equation*}
which satisfies \eqref{eq:RH}. The shocks are all admissible since they are increasing and $g$ is concave.

\emph{\textit{Repulsive case}}. The proof in this case is more involved, since the initial Dirac delta singularities are `smoothed' out immediately after $t=0$. On the other hand, in this case we have the following explicit formula for the $L^2$ gradient flow solution (see Theorem \ref{thm:explicit})
\begin{equation}\label{eq:solution_L2}
    X(s,t)=X_0(s) + t(2s-1).
\end{equation}
It is clear that $X(\cdot,t)$ has at most $N$ points of discontinuity. Let us set, as in \eqref{eq:rearrangement},
\begin{equation*}
    F(x,t)=\int_0^1 \chi_{(-\infty,x]}(X(s,t))ds.
\end{equation*}

Let us first prove that $F$ is a weak solution to the scalar conservation law \eqref{eq:CL}. We have to prove that
$F(x,t)$ solves, for all $\phi \in C_c^{\infty}([0,+\infty)\times \R)$,
\begin{equation}\label{eq:Fweaksolution}
\int_{\R}F_0\phi_0dx+\int_0^{+\infty}\int_{\R}\left[ F\phi_t+ (F^2-F)\phi_x \right] dxdt = 0.
\end{equation}

\begin{figure}[t]
\centering
\begin{tabular}{c  c }
\begin{tikzpicture}[scale=0.9]
\draw (0,0) node[anchor=south] {};
\draw (0,2.5) node[anchor=south] {$\sigma=1$};
\end{tikzpicture}
&
\begin{tikzpicture}[scale=0.9]
\draw[semithick, color=black] (-2,0) -- (-2,4);
\draw[semithick, color=black] (-2,0) -- (2,0);
\draw[semithick, color=black] (-2,-.1) -- (-2,0);
\draw (-2,-.1) node[anchor=north] {$0$};
\draw[semithick, color=black] (-1.5,-.1) -- (-1.5,.1);
\draw[semithick, color=black] (-1,-.1) -- (-1,.1);
\draw[semithick, color=black] (-.5,-.1) -- (-.5,.1);
\draw[semithick, color=black] (0,-.1) -- (0,.1);
\draw[semithick, color=black] (.5,-.1) -- (.5,.1);
\draw (-.5,-.1) node[anchor=north] {$1/2$};
\draw[semithick, color=black] (1,-.1) -- (1,.1);
\draw (1,-.1) node[anchor=north] {$1$};
\draw (1.9,0) node[anchor=north] {$s$};
\draw (-2,3.8) node[anchor=east] {$X_t(s)$};
\draw[semithick, color=black] (-2.1,1) -- (-1.9,1);
\draw[semithick, color=black] (-2.1,1.5) -- (-1.9,1.5);
\draw (-2,1.5) node[anchor=east] {$x_i$};
\draw[semithick, color=black] (-2.1,2) -- (-1.9,2);
\draw[semithick, color=black] (-2.1,3) -- (-1.9,3);
\draw[semithick, color=black] (-2,1) -- (-1,1);
\draw[->, semithick, color=black] (-1.5,.8) -- (-1.5,1.2);
\draw[semithick, color=black] (-1,1.5)  -- (-.5,1.5) ;
\draw[->, semithick, color=black] (-.75,1.4) -- (-.75,1.6);
\draw[semithick, color=black] (-.5,2) -- (.5,2) ;
\draw[->, semithick, color=black] (0,2.2) -- (0,1.8);
\draw[semithick, color=black] (.5,3)  -- (1,3);
\draw[->, semithick, color=black] (.75,2.8) -- (.75,2);
\draw[dashdotted, color=black] (-2,1.2) -- (-1,1.2);
\draw[dashdotted, color=black] (-1,1.6)  -- (-.5,1.6) ;
\draw[dashdotted, color=black] (-.5,1.8) -- (.5,1.8) ;
\draw[dashdotted, color=black] (.5,2.4)  -- (1,2.4);
\draw[dotted, color=black] (-2,1.4) -- (-1,1.4);
\draw[dotted, color=black] (-1,1.68)  -- (-.5,1.68) ;
\draw[dotted, color=black] (-.5,1.68) -- (.5,1.68) ;
\draw[dotted, color=black] (.5,1.9)  -- (1,1.9);
\draw[semithick, color=black] (7,0) -- (7,4);
\draw[semithick, color=black] (3,0) -- (11,0);
\draw[semithick, color=black] (7,-.1) -- (7,0);
\draw (10.9,0) node[anchor=north] {$x$};
\draw (7,3.8) node[anchor=east] {$F_t(x)$};
\draw (7,3) node[anchor=east] {$1$};
\draw[semithick, color=black] (7.1,3) -- (6.9,3);
\draw[semithick, color=black] (7.1,.5) -- (6.9,.5);
\draw[semithick, color=black] (7.1,1) -- (6.9,1);
\draw[semithick, color=black] (5,1)  -- (6,1);
\draw[semithick, color=black] (5,.1) -- (5,-.1);
\draw[semithick, color=black] (6,1.5)  -- (7.5,1.5);
\draw[semithick, color=black] (6,.1) -- (6,-.1);
\draw (6,0) node[anchor=north] {$x_i$};
\draw[semithick, color=black] (7.5,2.5)  -- (8.9,2.5);
\draw[semithick, color=black] (7.5,.1) -- (7.5,-.1);
\draw[semithick, color=black] (8.9,3)  -- (11,3);
\draw[semithick, color=black] (8.5,.1) -- (8.5,-.1);
\draw[semithick, color=black] (5,0) -- (5,1);
\draw[->, semithick, color=black] (4.9,.5) -- (5.7,.5);
\draw[semithick, color=black] (6,1) -- (6,1.5);
\draw[->, semithick, color=black] (5.9,1.25) -- (6.5,1.25);
\draw[semithick, color=black] (7.5,1.5) -- (7.5,2.5);
\draw[->, semithick, color=black] (7.6,2) -- (6.9,2);
\draw[semithick, color=black] (8.9,2.5) -- (8.9,3);
\draw[->, semithick, color=black] (9,2.75) -- (7.5,2.75);
\draw[dashdotted, color=black] (5.5,0)  -- (5.5,1);
\draw[dashdotted, color=black] (5.5,1)  -- (6.3,1);
\draw[dashdotted, color=black] (6.3,1)  -- (6.3,1.5);
\draw[dashdotted, color=black] (6.3,1.5)  -- (7.1,1.5);
\draw[dashdotted, color=black] (7.1,1.5)  -- (7.1,2.5);
\draw[dashdotted, color=black] (7.1,2.5)  -- (8,2.5);
\draw[dashdotted, color=black] (8,2.5)  -- (8,3);
\draw[dashdotted, color=black] (8,3)  -- (11,3);
\draw[dotted, color=black] (5.9,0)  -- (5.9,1);
\draw[dotted, color=black] (5.9,1)  -- (6.7,1);
\draw[dotted, color=black] (6.7,1)  -- (6.7,2.5);
\draw[dotted, color=black] (6.7,2.5)  -- (7.3,2.5);
\draw[dotted, color=black] (7.3,2.5)  -- (7.3,3);
\draw[dotted, color=black] (7.3,3)  -- (8,3);
\end{tikzpicture}
\\
\begin{tikzpicture}[scale=0.9]
\draw (0,0) node[anchor=south] {};
\draw (0,2.5) node[anchor=south] {$\sigma=-1$};
\end{tikzpicture}
&
\begin{tikzpicture}[scale=0.9]
\draw[semithick, color=black] (-2,0) -- (-2,4);
\draw[semithick, color=black] (-2,0) -- (2,0);
\draw[semithick, color=black] (-2,-.1) -- (-2,0);
\draw (-2,-.1) node[anchor=north] {$0$};
\draw[semithick, color=black] (-1.5,-.1) -- (-1.5,.1);
\draw[semithick, color=black] (-1,-.1) -- (-1,.1);
\draw[semithick, color=black] (-.5,-.1) -- (-.5,.1);
\draw[semithick, color=black] (0,-.1) -- (0,.1);
\draw[semithick, color=black] (.5,-.1) -- (.5,.1);
\draw (-.5,-.1) node[anchor=north] {$1/2$};
\draw[semithick, color=black] (1,-.1) -- (1,.1);
\draw (1,-.1) node[anchor=north] {$1$};
\draw (1.9,0) node[anchor=north] {$s$};
\draw (-2,3.8) node[anchor=east] {$X_t(s)$};
\draw[semithick, color=black] (-2.1,1) -- (-1.9,1);
\draw[semithick, color=black] (-2.1,1.5) -- (-1.9,1.5);
\draw (-2,1.5) node[anchor=east] {$x_i$};
\draw[semithick, color=black] (-2.1,2) -- (-1.9,2);
\draw[semithick, color=black] (-2.1,3) -- (-1.9,3);
\draw[semithick, color=black] (-2,1) -- (-1,1);
\draw[semithick, color=black] (-1,1.5)  -- (-.5,1.5) ;
\draw[semithick, color=black] (-.5,2) -- (.5,2) ;
\draw[semithick, color=black] (.5,3)  -- (1,3);
\draw[dashdotted, color=black] (-2,.7) -- (-1,.9);
\draw[dashdotted, color=black] (-1,1.4)  -- (-.5,1.5) ;
\draw[dashdotted, color=black] (-.5,2) -- (.5,2.2) ;
\draw[dashdotted, color=black] (.5,3.2)  -- (1,3.3);
\draw[dotted, color=black] (-2,.4) -- (-1,.8);
\draw[dotted, color=black] (-1,1.3)  -- (-.5,1.5) ;
\draw[dotted, color=black] (-.5,2) -- (.5,2.4) ;
\draw[dotted, color=black] (.5,3.4)  -- (1,3.6);
\draw[->, semithick, color=black] (-1.5,1.1) -- (-1.5,.6);
\draw[->, semithick, color=black] (-.75,1.6) -- (-.75,1.4);
\draw[->, semithick, color=black] (0,1.9) -- (0,2.2);
\draw[->, semithick, color=black] (.75,2.8) -- (.75,3.5);
\draw[semithick, color=black] (7,0) -- (7,4);
\draw[semithick, color=black] (3,0) -- (11,0);
\draw[semithick, color=black] (7,-.1) -- (7,0);
\draw (10.9,0) node[anchor=north] {$x$};
\draw (7,3.8) node[anchor=east] {$F_t(x)$};
\draw (7,3) node[anchor=east] {$1$};
\draw[semithick, color=black] (7.1,3) -- (6.9,3);
\draw[semithick, color=black] (7.1,.5) -- (6.9,.5);
\draw[semithick, color=black] (7.1,1) -- (6.9,1);
\draw[semithick, color=black] (5,1)  -- (6,1);
\draw[semithick, color=black] (5,.1) -- (5,-.1);
\draw[semithick, color=black] (6,1.5)  -- (7.5,1.5);
\draw[semithick, color=black] (6,.1) -- (6,-.1);
\draw (6,0) node[anchor=north] {$x_i$};
\draw[semithick, color=black] (7.5,2.5)  -- (8.9,2.5);
\draw[semithick, color=black] (7.5,.1) -- (7.5,-.1);
\draw[semithick, color=black] (8.9,3)  -- (11,3);
\draw[semithick, color=black] (8.5,.1) -- (8.5,-.1);
\draw[semithick, color=black] (5,0) -- (5,1);
\draw[semithick, color=black] (6,1) -- (6,1.5);
\draw[semithick, color=black] (7.5,1.5) -- (7.5,2.5);
\draw[semithick, color=black] (8.9,2.5) -- (8.9,3);
\draw[dashdotted, color=black] (4.4,0)  -- (4.8,1);
\draw[dashdotted, color=black] (4.8,1)  -- (5.8,1);
\draw[dashdotted, color=black] (5.8,1)  -- (6,1.5);
\draw[dashdotted, color=black] (7.5,1.5)  -- (7.9,2.5);
\draw[dashdotted, color=black] (8.9,2.5)  -- (9.3,2.5);
\draw[dashdotted, color=black] (9.5,3)  -- (9.3,2.5);
\draw[dotted, color=black] (3.8,0)  -- (4.6,1);
\draw[dotted, color=black] (4.6,1)  -- (5,1);
\draw[dotted, color=black] (5.6,1)  -- (6,1.5);
\draw[dotted, color=black] (7.5,1.5)  -- (8.3,2.5);
\draw[dotted, color=black] (8.9,2.5)  -- (9.7,2.5);
\draw[dotted, color=black] (9.7,2.5)  -- (10.1,3);
\draw[->, semithick, color=black] (5.1,.5) -- (4.2,.5);
\draw[->, semithick, color=black] (6.1,1.25) -- (5.8,1.25);
\draw[->, semithick, color=black] (7.4,2) -- (7.8,2);
\draw[->, semithick, color=black] (8.8,2.75) -- (9.7,2.75);
\end{tikzpicture}
\end{tabular}
\caption{Profile of $F_t$ and $X_t$ at $0=t_0 < t_1 < t_2$ ($t_0$ thick, $t_1$ dashdotted, $t_2$ dotted)}
\label{fig:mutft}
\end{figure}
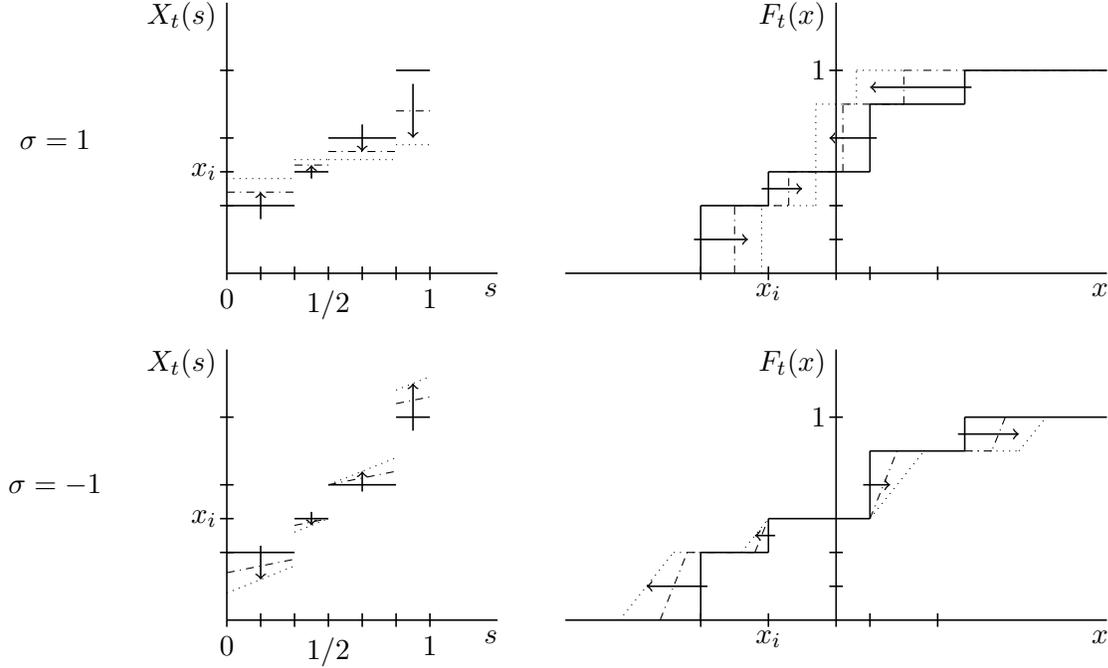

Let us set $J(t)=\{x\in \R\; |\;\; X(z,t)\neq x,\ \hbox{ for all }\ \ z\in [0,1]\}$, namely $J(t)$ is the complement of the image of $X(\cdot,t)$. Since $X(\cdot,t)$ has a finite number of jumps, $J(t)$ is the union of a finite number of disjoint open intervals. It is easily seen that $F$ is constant along each connected component of $J(t)$. Now, let $x\in (\R\setminus J(t))^\circ$. We have
\begin{equation}\label{eq:der-F}
\partial_x F(x,t)=
\left( \partial_s X \right)^{-1}(F(x,t))= \frac{1}{2t},
\end{equation}
because $X(\cdot,t)$ is monotonic increasing on a small neighborhood of $x$. Therefore, $\mu(\cdot,t):=\frac{\partial F}{\partial x}(\cdot,t)$ is absolutely continuous for all $t>0$ on each component of $(\R\setminus J(t))^\circ$. Moreover, in a small neighborhood of $x\in (\R\setminus J(t))^\circ$ we have
\begin{equation}\label{eq:der-identity}
0=\frac{d}{dt}X(F(x,t),t)=\partial_t X(F(x,t),t)+\partial_s X(F(x,t),t) \partial_t F(x,t).
\end{equation}
Finally, let us notice that $F(\cdot,t)$ is absolutely continuous on $\R$, and therefore it is differentiable w.r.t. $x$ almost everywhere, with $\mu(\cdot,t):=\partial_x F(\cdot,t)$ being a probability measure (see Figure~\ref{fig:mutft} for the difference between attractive and repulsive case).

Then, for each $\phi \in C_c^{\infty}$ we have
\begin{equation*}
\begin{split}
& \int_0^{+\infty}\int_{\R}(\partial_x \phi) (F^2-F)dxdt=-\int_0^{+\infty}\int_{\R}\phi(2F-1)F_xdxdt=-\int_0^{+\infty}\int_{\R}\phi(2F-1)d\mu dt  \\
& = -\int_0^{+\infty}\int_0^1\phi (X(s,t))(2s-1)dsdt = -\int_0^{+\infty}\int_0^1\phi(X(s,t),t)\partial_t Xdsdt\\
& = -\int_0^{+\infty}\int_{J(t)^c}\phi(x,t) \left.\partial_t X\right|_{s=F(x,t)}\partial_xFdxdt,
\end{split}
\end{equation*}
where we use \eqref{eq:solution_L2}. Now, using \eqref{eq:der-F} and \eqref{eq:der-identity}, we obtain
\begin{equation*}
\begin{split}
& -\int_0^{+\infty}\int_{J(t)^c}\phi(x,t) \left.\partial_t X\right|_{s=F(x,t)}\partial_xFdxdt \\
& = \int_0^{+\infty}\int_{J(t)^c}\phi(x,t) \partial_t F \left. \partial_s X \right|_{s=F(x,t)}  \partial_x F dx dt = \int_0^{+\infty}\int_{J(t)^c}\phi(x,t) \partial_t F dx dt \\
& = \int_0^{+\infty}\int_{\R}\phi(x,t) \partial_t F dx dt = -\int_{\R}\phi_0(x)F_0(x)dx -\int_0^{+\infty}\int_{\R} (\partial_t\phi) F dxdt,
\end{split}
\end{equation*}
which proves the assertion \eqref{eq:Fweaksolution}.

We next prove that $F(x,t)$ satisfies the Oleinik condition \eqref{Oleinik}. Given any $0 \leq s_1 < s_2 \leq 1$ we have that
\begin{equation} \label{eq:satisfyoleinik}
X(s_2,t)-X(s_1,t)=X(s_2,0)-X(s_1,0)+2t(s_2-s_1) \implies \; (s_2-s_1) \leq \frac{X(s_2,t)-X(s_1,t)}{2t}.
\end{equation}
From the definition of $F$ in terms of his pseudo-inverse \eqref{eq:rearrangement} we obtain that
\begin{equation*}
F(x+z,t)-F(x,t)=\int_0^1 \mathcal{X}_{(x,x+z]}(X(s,t))ds= |\{ x< X(s,t) \leq x+z \}| =s_2^*-s_1^*,
\end{equation*}
where
\begin{equation*}
\left. \begin{cases}
s_2^*= \sup \left\{s| X(s,t) \in (x,x+z] \right\}\\
s_1^*= \inf \left\{s| X(s,t) \in (x,x+z] \right\}
\end{cases} \right.  .
\end{equation*}
Using \eqref{eq:satisfyoleinik} we deduce
\begin{equation*}
F(x+z,t)-F(x,t)=s_2^*-s_1^* \leq \frac{X(s_2^*,t)-X(s_1^*,t)}{2t} \leq \frac{x+z - x}{2t} = \frac{z}{2t}.
\end{equation*}
and then the Oleinik condition \eqref{Oleinik} is satisfied. \\
We now prove that $\mu(x,t)=\frac{\partial F}{\partial x}$ is the solution of \eqref{eq:main} satisfying \eqref{GradFlow}. Let us first see that $\mu(x,t)$ is a weak measure solution of the continuity equation with the velocity field $v(x,t)=2F-1 $. With $\Phi(x,t) \in C_c^{\infty}([0,+\infty]\times \R)$, by direct integration by parts we obtain
\begin{align*}
 \int_0^{+\infty}\!\!\int_{\R} &\left(  \partial_x\Phi \, (2F-1)+\partial_t \Phi \right)d\mu dt+\int_{\R}\Phi_0 d\mu_0 \\
=&\int_0^{+\infty}\int_{\R} \left(\partial_x\Phi \, (2F-1)+\partial_t\Phi \right) \partial_x F dx dt+\int_{\R}\Phi_0 \partial_xF_0dx  \\
=& \int_0^{+\infty}\int_{\R} \partial_x\Phi \, \partial_x(F^2-F) dxdt + \int_0^{+\infty}\int_{\R} \partial_t\Phi \, \partial_x F dx dt+\int_{\R}\Phi_0 \partial_x F_0 dx  \\
=&-\int_0^{+\infty}\int_{\R} \partial_x^2\Phi \, (F^2 -F)dxdt  -\int_0^{+\infty}\int_{\R}\partial_t \partial_x  \Phi \, F dx dt-\int_{\R}\partial_x\Phi_0 \, F_0dx.
\end{align*}
Now, choosing $\phi=-\partial_x \Phi$, \eqref{eq:Fweaksolution} implies that the previous equation is equal to
\begin{equation*}
\int_0^{+\infty}\int_{\R} \partial_x\phi \, (F^2 -F)dxdt+\int_0^{+\infty}\int_{\R}\partial_t\phi \, F dx dt+\int_{\R}\phi_{0} F_0dx = 0,
\end{equation*}
and so $\mu(x,t)$ is a weak solution.
The second condition \eqref{GradFlow} comes straightforwardly:
\begin{equation*}
v(x,t)=(2F-1)=\int_{x \neq y}\sign(x-y)d\mu(y,t)=-\partial^0 \W [\mu].
\end{equation*}

\textsc{Step 2. General initial measure}

Let $\mu_0\in \PP$ be a general initial condition. Define $F_0(x):=\mu_0((-\infty,x])$, and let $X_0$ be the pseudo-inverse of $F_0$. We then denote by $\mu_t$ the unique Wasserstein gradient flow solution to \eqref{eq:main} with initial condition $\mu_0$, by $\widetilde{F}(\cdot,t)$ the unique entropy solution to \eqref{eq:CL}, and by $\widetilde{X}_t$ the unique $L^2$ gradient flow solution of \eqref{eq:GF_L2}. As usual we set $F(x,t)=\mu_t((-\infty,x])$ and by $X_t\in L^2({0,1})$ the pseudo-inverse of $F(\cdot,t)$. We need to prove that $F=\widetilde{F}$ and $X=\widetilde{X}$. Let $\mu_0^N\in \PP$ be a linear combination of Dirac Delta as in Step 1, such that $d_W(\mu_0^N,\mu_0)\rightarrow 0$ as $N\rightarrow +\infty$. Let us recall the definition of $1$-Wasserstein distance between $\nu_1,\nu_2\in \PP$
\begin{equation}
\label{def:W1}
  d_{W,1}(\nu_1,\nu_2)=\inf\left\{\int_{\R\times \R}|x-y|d\bm{\gamma}(x,y),\;\;\bm{\gamma} \in \Gamma(\nu_1,\nu_2)\right\} = \|F_1 - F_2\|_{L^1(\R)}\,,
\end{equation}
with $F_i(x)=\nu_i((-\infty,x])$.
For $\nu_1,\nu_2\in \PP$ it is immediately seen that $d_{W,1}(\nu_1,\nu_2)\leq  d_{W}(\nu_1,\nu_2)$.
Let $F_0^N(x)=\mu_0^N((-\infty,x])$ and $F^N(\cdot,t)$ be the unique entropy solution to \eqref{eq:CL} with initial condition $F_0^N$. Let $X_0^N$ be the pseudo-inverse of $F_0^N$ and let $X_t^N$ be the unique $L^2$ gradient flow solution to \eqref{eq:GF_L2} with initial condition $X^N_0$. Due to \eqref{eq:contraction_wass}, we have for all $t\geq 0$
\begin{align*}
  & \|F(t)-F^N(t)\|_{L^1} = d_{W,1}(\mu(t),\mu^N(t)) \leq d_{W}(\mu(t),\mu^N(t))\leq d_{W}(\mu_0,\mu^N_0) \ \rightarrow 0 \;\; \hbox{ as }\;\; N\rightarrow +\infty.
\end{align*}
Moreover, from the $L^1$ contraction \eqref{eq:contraction_entropy} in theorem \ref{thm:entropy} we get
\begin{align*}
   & \|\widetilde{F}(t)-F^N(t)\|_{L^1}\leq \|F_0-F^N_0\|_{L^1} \ \rightarrow 0 \;\; \hbox{ as }\;\; N\rightarrow +\infty.
\end{align*}
By uniqueness of the limit, $\widetilde{F}\equiv F$ for all $t\geq 0$. Similarly, from the definition of $d_W$ \eqref{eq:distance_rappr} and the $d_W$ contraction \eqref{eq:contraction_wass} in theorem \ref{thm:dwcontraction}, we get
\begin{equation*}
  \|X_t - X_t^N\|_{L^2} = d_W (\mu_t,\mu^N_t)\leq  d_W (\mu_0,\mu^N_0)\ \rightarrow 0 \;\; \hbox{ as }\;\; N\rightarrow +\infty.
\end{equation*}
Finally, from \eqref{eq:contraction_L2} we get
\begin{equation*}
  \|\widetilde{X}_t-X_t^N\|_{L^2}\leq  \|X_0-X_0^N\|_{L^2} = d_W (\mu_0,\mu^N_0) \ \rightarrow 0 \;\; \hbox{ as }\;\; N\rightarrow +\infty,
\end{equation*}
and the assertion follows.
\endproof

\section{Particle approximation}\label{sec:particle}

A clear distinction between the attractive case $\sigma =1$ and the repulsive case $\sigma=-1$ is that the former case allows for \emph{atomic measure} solutions as a special case of gradient flow solutions, whereas this is not possible in the latter case. More precisely, in the attractive case, assuming
\begin{equation}\label{eq:initial_deltas}
     \mu_0 = \sum_{j=1}^N m_j\delta_{x_j},
\end{equation}
as in Step 1 of Theorem \ref{thm:equivalence}, if the vector $(x_j(t))_{j=1}^N$ is the (unique) solution to the particle system
\begin{equation*}
    \dot{x}_j(t)= - \sum_{k=1}^N m_k \sign(x_j(t)-x_k(t)),  \qquad \qquad (\sign(0) =0),
\end{equation*}
then, as we proved in Theorem~\ref{thm:equivalence}, the empirical measure $\mu_t=\sum_{j=1}^N m_j\delta_{x_j(t)}$ is the unique gradient flow of $\W$ with $\sigma=1$ with initial condition $\mu_0$. By the stability property \eqref{eq:contraction_wass} then allows any gradient flow solution $\mu_t$ to be approximated by the empirical measure of a finite number of particles, uniformly in time. Note that the approximating empirical measures are \emph{exact solutions} of the same problem.

On the other hand, in the repulsive case $\sigma=-1$ the proof of Theorem \ref{thm:equivalence} shows that the unique gradient flow $\mu_t$ of $\W$ with initial condition $\mu_0$ as in \eqref{eq:initial_deltas} is the  $x$-derivative of a continuous piecewise linear function $F(\cdot,t)$ consisting of $N$ rarefaction waves. Hence $\mu_t$ is absolutely continuous with respect to the Lebesgue measure. Therefore, in the repulsive case the approximation of an arbitrary gradient flow solution by a finite number of moving deltas is not as simple as in the attractive case. In the next theorem we provide a solution to such problem, which was recently addressed also in \cite{CCH13} for a class of singular interaction potentials in many space dimensions. In our specific case, it turns out that the particle approximation for \eqref{eq:main} is equivalent to the convergence of the so called wave front tracking scheme for the scalar conservation law \eqref{eq:CL}, see \cite{Daf72,dip76,Bre92}.

In order to state the result we introduce some notation. For a given initial probability measure $\mu_0 \in \mathcal{P}(\R)$ and  a fixed positive integer $N$, we define inductively the finite sequence $\{X^N_j\}_{j=1}^{N}$ as
\begin{equation}\label{eq:particle_initial}
  \begin{cases}
  X^N_1 = \inf\left\{x\in \R\ :\; \mu_0((-\infty,x))>\frac{1}{2N}\right\}  \\
  X^N_{j+1}=\inf\left\{x\in \R\ :\; \mu_0([X_j^N,x))>\frac{1}{N}\right\}
  \end{cases} j=1,\ldots, N-1.
\end{equation}
Roughly speaking, we have divided the support of $\mu_0$ into $N$ sets on which $\mu_0$ has equal mass $1/N$, and chosen the position $X_j^N$ to be an intermediate point of that interval. Such a construction could be much easier in the case of $\mu_0$ with bounded support (e.g. by assigning the position of each particle on the edge of each mass portion), but we choose this construction to include initial data with unbounded support. Next, we define the empirical measure
\begin{equation*}
  \mu_0^N=\sum_{j=1}^{N} \frac{1}{N}\delta_{X^N_j}.
\end{equation*}
We define the cumulative distribution function of $\mu^N_0$ as
\begin{equation*}
  F_0^N=\sum_{j=1}^{N} \frac{1}{N}\chi_{[X_j^N,+\infty)}.
\end{equation*}
We now introduce the approximated flux
\begin{align*}
 g_N(F)& \ =\sum_{j=1}^{N}N\left[g\left(\frac{j}{N}\right)-g\left(\frac{j-1}{N}\right)\right]
  \left(F-\frac{j-1}{N}\right)\chi_{[\frac{j}{N},\frac{j-1}{N})}(F) .
\end{align*}

\begin{figure}[h]
\centering
\begin{tikzpicture}
\draw[dashed,domain=-2:2] plot (\x,{7*(0.25*\x+.5)*(0.25*\x+.5)-1.75*\x-1});
\draw[semithick, color=black] (-2,1) -- (-2,4);
\draw[semithick, color=black] (-2.3,2.5) -- (2.3,2.5);
\draw[semithick, color=black] (-2,2.5) -- (-1,{7*(0.25-.5)*(0.25-.5)+1.75-1});
\draw[semithick, color=black] (-1,2.4) -- (-1,2.6);
\draw[semithick, color=black] (-1,{7*(0.25-.5)*(0.25-.5)+1.75-1}) -- (0,.75);
\draw[semithick, color=black] (0,2.4) -- (0,2.6);
\draw[semithick, color=black] (0,.75) -- (1,{7*(0.25+.5)*(0.25+.5)-1.75-1});
\draw[semithick, color=black] (1,2.4) -- (1,2.6);
\draw[semithick, color=black] (1,{7*(0.25+.5)*(0.25+.5)-1.75-1}) -- (2,2.5);
\draw[semithick, color=black] (2,2.4) -- (2,2.6);
\draw (-2,2.2) node[anchor=east] {$0$};
\draw (-2,2.2) node[anchor=east] {$0$};
\draw (2,2.4) node[anchor=north] {$1$};
\draw (2.3,2.5) node[anchor=north] {$F$};
\draw (-2,3.8) node[anchor=east] {$g^N(F)$};
\draw[semithick, color=black] (7,0.9) -- (7,4);
\draw[semithick, color=black] (2.5,1) -- (11,1);
\draw (10.9,1) node[anchor=north] {$x$};
\draw (7,3.8) node[anchor=east] {$F^N_0(x)$};
\draw (7,3) node[anchor=east] {$1$};
\draw[semithick, color=black] (7.1,3) -- (6.9,3);
\draw[semithick, color=black] (7.1,2) -- (6.9,2);
\draw[semithick, color=black] (7.1,1.5) -- (6.9,1.5);
\draw[semithick, color=black] (7.1,2.5) -- (6.9,2.5);
\draw[dashed, color=black] (4,1) -- (9,3);
\draw[semithick, color=black] (4.625,1.1) -- (4.625,0.9);
\draw[semithick, color=black] (4.625,1.5)  -- (5.775,1.5) circle (2pt);
\fill[color=black] (4.625,1.5) circle (0.3ex);
\draw[semithick, color=black] (5.775,1.1) -- (5.775,0.9);
\draw[semithick, color=black] (5.775,2)  -- (7.025,2) circle (2pt);
\fill[color=black] (5.775,2) circle (0.3ex);
\draw[semithick, color=black] (7.025,1.1) -- (7.025,0.9);
\draw[semithick, color=black] (7.025,2.5)  -- (8.275,2.5) circle (2pt);
\fill[color=black] (7.025,2.5) circle (0.3ex);
\draw[semithick, color=black] (8.275,1.1) -- (8.275,0.9);
\draw[semithick, color=black] (8.275,3)  -- (9,3);
\fill[color=black] (8.275,3) circle (0.3ex);
\draw[semithick, color=black] (9,3)  -- (11,3);
\end{tikzpicture}
\caption{A pictorial view of the approximated problem with $N=4$}
\label{fig55}
\end{figure}
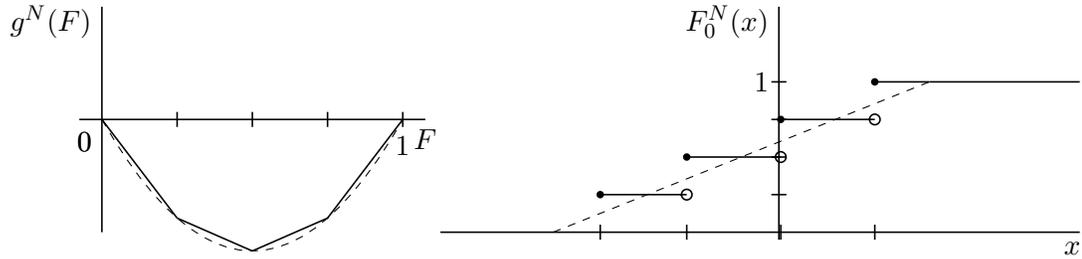
Figure \ref{fig55} illustrates the construction of $g^N$. Notice that $g^N$ is piecewise linear and convex on $[0,1]$. We now define the approximating distribution $F^N(x,t)$ as the unique $L^\infty$ solution to
\begin{equation}\label{eq:WFT}
  \partial_t F^N + \partial_x g^N(F^N) = 0,
\end{equation}
with initial condition $F^N_0$. Let $\mu^N$ be the $x$-derivative in the sense of distributions of the solution $F^N(x,t)$ of \eqref{eq:WFT}.

The solution $F^N$ to \eqref{eq:WFT} consists of exactly $N$ shocks, with constant velocities defined by the Rankine-Hugoniot condition
\begin{equation}\label{eq:RH1}
\dot x_j(t) = \lambda_j:= \frac{g\left(\frac{j}{N}\right)-g\left(\frac{j-1}{N}\right)}{\frac{1}{N}} = \frac{(j-1)}{N} - \frac{(N-j)}{N},
\end{equation}
and with initial positions $x_j(0)=X^N_j$, as was first observed in~\cite[Section 6]{Bre92}. We have then the explicit formula for the shock curves
\begin{equation*}
  x_j(t)=X^N_j + \left[\frac{(j-1)}{N} - \frac{(N-j)}{N}\right] t,
\end{equation*}
and the explicit formula for the solution $F^N$ is given by
\begin{equation*}
F^N(x,t)= \sum_{j=1}^{N-1} \chi_{[x_j(t),x_{j+1}(t))}(x) \frac{j}{N} + \chi_{[x_N(t),+\infty)}(x).
\end{equation*}
In the evolution of $x_j(t)$ the shocks do not cross each other, since $\lambda_j < \lambda_{k}$ if $j<k$. Also note that all the shocks have the same size in the jump, namely $1/N$. This is consistent with the fact that no shocks will appear in the continuum limit, as the flux is convex. Moreover, the formula \eqref{eq:RH1} shows that each discontinuity $x_j$ is driven by a positive drift $(j-1)/N$, which can be interpreted as a repulsive force against the $j-1$ particles on its left, and a negative drift $(N-j)/N$, which accounts for a repulsive force against the $N-j$ particles on its right. Note that we introduced an ordering between the particles; consider for example the situation where the starting point is a single shock. At time $t = 0$ there is no notion of particles on the left/right but still the evolution, according to the Rankine-Hugoniot condition, prescribes a velocity $\lambda_j$ to the $j$-th particle.

In the next Theorem we prove that the empirical measure $\mu^N(t)$ converges in the $2$-Wasserstein distance to the solution $\mu$ to the repulsive gradient flow. In the landscape of conservation laws, this is equivalent to prove that $F^N(t)$ converges in some sense to the cumulative distribution $F(x,t)$ of $\mu(t)$ (convergence in $L^1$ of $F^N$ means convergence in the $1$-Wasserstein distance of $\mu^N$ to $\mu$). One way to perform this task could then be using the same strategy of \cite{Daf72,dip76,Bre92}, which relies on providing $BV$ estimates on $F^N$. However, in our case we have explicit formulas for $\mu^N(t)$ and $\mu(t)$, so the convergence can be checked directly.

\begin{theorem}[Particle approximation in the repulsive case]\label{thm:particle}
Let $\mu_0 \in \mathcal{P}_2(\R)$, and let $\mu(x,t)$ be the unique gradient flow solution of $\W$ with $\sigma=-1$ with initial datum $\mu_0$.
For each $N$, let $\mu^N$ be the empirical measure
\begin{equation*}
  \mu^N(t)= \sum_{j=1}^{N} \frac{1}{N}\delta_{x^N_j(t)},
\end{equation*}
with $x^N_j$ satisfying
\begin{equation}\label{eq:particle_ODE}
\dot{x}^N_j(t)=\frac1N \sum_{k\neq j} \sign(x_j^N(t)-x_k^N(t)) = \frac{2j-1-N}{N},  \qquad x_j^N(0) = X^N_j, \qquad j=1,\cdots,N,
\end{equation}
where $X^N_j$ is defined in \eqref{eq:particle_initial}. Then, for all $t\geq 0$, we have
\begin{equation*}
\lim_{N\to \infty}d_W(\mu^N(t),\mu(t))=0.
\end{equation*}
\end{theorem}

\proof

From a direct computation with pseudo inverses, we can easily check that the pseudo inverse variable $X^N$ related to the empirical measure $\mu^N$ can be written as
\begin{align*}
& X^N(s,t)=X_0^N + f^N(s)t,
\end{align*}
with $f^N(s)=\frac{2j-1-N}{N}$ for $s \in [\frac{j-1}{N},\frac{j}{N})$. Moreover, we recall that the pseudo inverse $X$ related to $\mu$ can be written as
\begin{align*}
    & X(s,t)=X_0+(2s-1)t.
\end{align*}
Therefore, we obtain that
\begin{align*}
d_W^2(&\,\mu(t),\mu^N(t)) = \| X(t)-X^N(t) \|^2_{L^2} = \int_0^1 \left( X_0+(2s-1)t-X^N_0-f^N(s)t \right)^2 ds \\
& \,= \int_0^1 \left( X_0-X^N_0 \right)^2 ds + 2t\int_0^1 \left( X_0-X^N_0 \right) \left( 2s -1 -f^N(s) \right) ds + \int_0^1 \left(2s -1 -f^N(s)\right)^2 ds.
\end{align*}
Combining the previous equality with the bound $|2s -1 -f^N(s)| \leq \frac{1}{N}$ and with the convergence of $X_0^N$ to $X_0$ we conclude the proof.
\endproof

\begin{remark}
In our construction we chose a specific way to approximate the initial datum via a combination of deltas, namely by placing the particle at the mid-mass-point (the `mass median') in each interval. It can be easily checked that such construction is not necessary, and the convergence result works for more general approximation procedures for the initial datum.
\end{remark}

\section{The characterization of the sub-differential of $\W$}\label{sec:subdiff}

Here we analyse the sub-differential of the functional $\W$ in the repulsive case, namely with $\sigma=-1$. Let us remark here that this task is completely solved in the attractive case in view of the results in \cite{CDFLS}. We will extensively use that
\begin{equation}\label{funcdist}
\W[\mu]=\frac{1}{2}\int_{\R \times \R}|x-y|d\mu(x)d\mu(y)=\int_0^1 X_{\mu}(s)(2s-1)ds.
\end{equation}
From Proposition~\ref{0geodesicconvex} we have that, if we deal with a measure $\mu$ with no atoms, then the sub-differential is characterized as follows
\begin{equation}\label{eq:ridefinition}
 \partial^0\W(\mu) = - \int_{x \neq y}\sign(x-y) d\mu(y)=: k(x).
\end{equation}
In the case that $\mu$ has concentrated mass, then the sub-differential is empty, as proven in the following
\begin{theorem}
\label{th:charsubdiff}
Let $\mu \in \PP$ and $W(x)=-|x|$. If there exists $\overline{x} \in \R$ such that $\mu(\{ \overline{x} \}) > 0$ then $\partial \mathcal{W} (\mu) = \emptyset$. Conversely, if $\mu(\{x\})=0$ for all $x$, then $\partial^0 \W$, the element of minimal norm, is
\[
\partial^0\W(\mu) = - \int_{x \neq y}\sign(x-y) \,d\mu(y).
\]
\proof
The proof of the second statement can be found in \cite[Proposition~4.3.3]{BS11}.
We now prove the first statement. {\def\sigma{a} 
Assume that there exists $\overline{x}\in\R$ such that $\mu(\{ \overline{x} \})= \sigma > 0$. Then there exist $0 \leq r_1 < r_2 \leq 1$ such that
$
r_2=r_1+\sigma, \; X_{\mu}(s)\equiv \overline{x} \; \text{ for every }s\in [r_1,r_2].
$
We take $r_2$ to be maximal and $r_1$ minimal, i.e.,
\begin{equation}
\label{def:r2-maximal}
\text{for all }\delta>0, \quad X_\mu(r_2+\delta) > X_\mu(r_2-) \quad \mbox{ and }
\quad X_\mu(r_1-\delta) < X_\mu(r_1+).
\end{equation}

Assume that $\partial \W(\mu)$ is not empty; let $k\in L^2(\mu)$ be any element of $\partial \W$. For every measure $\nu \in \PP$, we have
\begin{align}\label{tech}
\W[\nu]-\W[\mu]&- \int_{\R \times \R}k(x) (y-x)\,d\bm{\gamma}(x,y) =\int_0^1\left( X_{\nu}(s)-X_{\mu}(s)\right)\left(1-2s-k(X_{\mu}(s)) \right) \,ds \,,
\end{align}
since $\bm{\gamma}$ is the optimal plan as in \eqref{eq:distance_rappr} taking into account \eqref{funcdist}.

We will arrive to a contradiction by constructing sequences of $\nu_\e$, converging to $\mu$ in $d_W$ as $\e\to0$, leading to conditions on $k$ that cannot be satisfied.

Given $\e>0$ and $0 < \eta < \sigma$ we define
\[
\delta_{\e}:= \inf \left\{ \theta \geq 0 | X_{\mu}(r_2+\theta) \geq \overline{x}+\e \right\}.
\]
It follows from~\eqref{def:r2-maximal} that $\delta_{\e} \to 0$ as $\e \to 0$ with $0\leq \delta_\e <1$.
Define $\nu$ by setting $X_\nu$ as follows:
\[
X_\nu(s) := \begin{cases}
\overline{x}+\e & \text{if } s\in [r_2-\eta,r_2+\delta_{\e}]\\
X_\mu(s) & \text{otherwise}.
\end{cases}
\]
By the definition of $\delta_\e$, this $X_\nu$ is increasing and therefore $\nu$ is well-defined. Although $\nu$ depends on $\e$, we do not indicate this to alleviate the notation. We calculate
\[
d_W^2(\mu,\nu) = \int_0^1 |X_\mu(s)-X_\nu(s)|^2 \, ds =
\int_{r_2-\eta}^{r_2+\delta_\e} |X_\mu(s)-X_\nu(s)|^2 \, ds ,
\]
implying that $d_W^2(\mu,\nu)\in [\e^2\eta,\e^2(\eta+1)]$.
Therefore, from \eqref{tech} we deduce that
\begin{align*}
\W[\nu]-\W[\mu]&- \int_{\R \times \R}k(x) (y-x)\,d\bm{\gamma}(x,y) \\
&=\int_{r_2-\eta}^{r_2}\e\left( 1 - 2s - k(\overline{x}) \right) ds + \int_{r_2}^{r_2+\delta_{\e}}\left(\overline{x}+\e-X_{\mu}(s)\right)\left( 1 - 2s - k(X_{\mu}(s)) \right) ds.
\end{align*}
We estimate the last term by
\begin{align*}
 \int_{r_2}^{r_2+\delta_{\e}}| \overline{x}+ \e &-X_{\mu}(s) | \; | 1 - 2s - k(X_{\mu}(s))|\,ds \leq \e \int_{r_2}^{r_2+\delta_{\e}}| 1 - 2s - k(X_{\mu}(s))|\,ds \\
 &\leq  \e \delta_{\e} +  \e\int_{r_2}^{r_2+\delta_{\e}}  |k(X_{\mu}(s))|ds
   \leq  \e \Bigl(\delta_{\e} + \sqrt{\delta_\e}\| k \|_{L^2(\mu)} \Bigr) .
\end{align*}
In order for $k$ to satisfy~\eqref{ineq:subdiffineq}, it is therefore necessary that
\[
\int_{r_2-\eta}^{r_2}\left( 1 - 2s - k(\overline{x}) \right) ds =  \eta (1 + \eta -2r_2 - k(\overline{x}))\geq \delta_{\e} + \sqrt{\delta_\e}\| k \|_{L^2(\mu)} ,
\]
which implies $k(\overline{x})\leq 1-2r_2+\eta$.
Note that this inequality applies for each choice of $0<\eta<\sigma$, and therefore we find that
$
k(\overline{x}) \leq 1-2r_2.
$
By repeating the argument for a similar interval $[r_1-\tilde \delta_{\e},r_1+\eta]$ we find a similar bound on $k(\overline x)$ from below. Together these inequalities read
$
1-2r_1 \leq k(\overline x) \leq 1- 2r_2.
$
Since $r_2> r_1$, it is impossible to satisfy both inequalities simultaneously, and we therefore find a contradiction.
} 
\endproof
\end{theorem}

\subsection{Extended sub-differential}

Let us recall the notion of extended subdifferential, more details can be found in \cite[Chapter 10]{AGS}. We define the set of optimal $3$-plans $\Gamma_0(\bm{\mu},\mu_3)$, for $\bm{\mu}\in \mathcal{P}_2(\R\times \R)$ and $\mu_3\in \PP$, as follows: $\bm{\gamma} \in \Gamma_0(\bm{\mu},\mu_3)$ if and only if $(\pi^{1,3})_\sharp\bm{\gamma} \in \Gamma_0((\pi_1)_\sharp \bm{\mu},\mu_3)$. Here $\pi^{1,3}$ is the projection of $\R\times\R\times \R$ onto the first and third components.

\begin{definition}[Extended Fr\'{e}chet sub-differential]\label{def:frechet2}
Let $\phi:\PP\rightarrow (-\infty,+\infty]$ be proper and lower semi-continuous, and let $\mu\in D(\phi)$. We say that $\bm{\gamma} \in \mathcal{P}_2(\R\times \R)$ belongs to the \emph{extended Frech\'{e}t sub-differential} $\bm{\partial}\phi(\mu)$ if $(\pi_1)_\sharp \bm{\gamma} =\mu$ and
\begin{equation}
\label{ineq:subdiffineq}
    \phi(\mutil)-\phi(\mu)\geq \inf_{\bm{\mu}\in \Gamma_0(\bm{\gamma},\mutil)}\int_{\R\times\R\times\R} x_2(x_3-x_1) d\bm{\mu}(x_1,x_2,x_3)+ o(d_W(\mu,\mutil)).
\end{equation}
\end{definition}

Assume $(i\otimes T)_{\#} \mu \in \Gamma_0(\mu,\mutil)$, with $\mutil = T_{\#} \mu$, which means that $T$ is an optimal map between $\mu$ and $\mutil$. Then, for each element $k$ of the Fr\'echet sub-differential we can construct an element of the extended sub-differential through the formula $\bm{\mu}=(i \otimes k \otimes T)_{\#}\mu$. In ~\cite[Chapter~10]{AGS} the authors prove the existence of an element of the extended sub-differential for a wide class of functionals called `regular' functionals. Such element may be not an element of the standard sub-differential, which may indeed be empty as in the present case.

Let $\bm{\gamma} \in \bm{\partial}\phi (\mu)$ be a plan $\bm{\gamma}\in \mathcal{P}_{2}(\R \times \R)$ such that $\pi_{\#}^1\bm{\gamma}=\mu$ and
\begin{equation*}
\W(\nu)-\W(\mu) \geq \inf_{\bm{\mu}\in \Gamma_0(\bm{\gamma},\nu)}\int_{\R^3}  x_2(x_3-x_1) d\bm{\mu} + o\bigl( d_W(\mu,\nu) \bigl).
\end{equation*}
For $\bm{\gamma} \in \mathcal{P}_{2}(\R \times \R)$ we need to define
\[
|\bm{\gamma}|^2_j := \int_{\R^2}|x_j|^2d\bm{\gamma}(x_1,x_2),  \qquad j=1,2.
\]
It is important to notice that, when $\mu$ is absolutely continuous w.r.t. the Lebesgue measure, then $k \in L^2(\mu)$ belongs to the Fr\'echet sub-differential $\partial \W(\mu)$ if and only if
\begin{equation*}
\bm{\gamma}=(i \otimes k)_{\#}\mu\in \bm{\partial}\W(\mu),
\end{equation*}
in fact $\Gamma_0(\bm{\gamma},\nu)$ is known in this case and it contains the unique element  $\bm{\mu}= (i \otimes k\otimes t_{\mu}^{\nu})_{\#}\mu$, for a more detailed discussion we refer to \cite[Chapter~10]{AGS}.

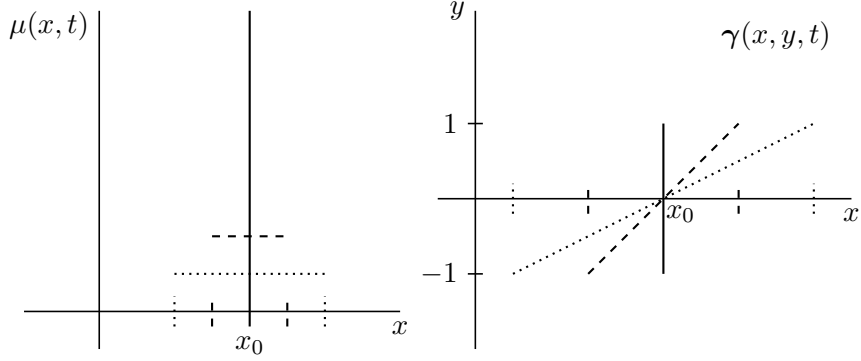
\begin{figure}[t]
\centering
\begin{tikzpicture}
\draw[semithick, color=black] (-2,-.5) -- (-2,4);
\draw[semithick, color=black] (-3,0) -- (2,0);
\draw (2,0) node[anchor=north] {$x$};
\draw (-2,3.8) node[anchor=east] {$\mu(x,t)$};
\draw[thick, color=black] (0,-.2) -- (0,4);
\draw (0,-.2) node[anchor=north] {$x_0$};
\draw[dashed, thick, color=black] (-.5,1) -- (.5,1);
\draw[dashed, thick, color=black] (-.5,-.2) -- (-.5,.2);
\draw[dashed, thick, color=black] (.5,-.2) -- (.5,.2);
\draw[dotted, thick, color=black] (-1,.5) -- (1,.5);
\draw[dotted, thick, color=black] (-1,-.2) -- (-1,.2);
\draw[dotted, thick, color=black] (1,-.2) -- (1,.2);
\draw (7,4) node[anchor=north] {$\bm{\gamma}(x,y,t)$};
\draw (8,1.5) node[anchor=north] {$x$};
\draw (3,4) node[anchor=east] {$y$};
\draw[semithick, color=black] (3,-.5) -- (3,4);
\draw[semithick, color=black] (2.5,1.5) -- (8,1.5);
\draw[semithick, color=black] (2.9,.5) -- (3.1,.5);
\draw (2.9,2.5) node[anchor=east] {$1$};
\draw[semithick, color=black] (2.9,2.5) -- (3.1,2.5);
\draw (2.9,.5) node[anchor=east] {$-1$};
\draw[thick, color=black] (5.5,.5) -- (5.5,2.5);
\draw (5.4,1.3) node[anchor=west] {$x_0$};
\draw[dashed, thick, color=black] (4.5,.5) -- (6.5,2.5);
\draw[dashed, thick, color=black] (4.5,1.3) -- (4.5,1.7);
\draw[dashed, thick, color=black] (6.5,1.3) -- (6.5,1.7);
\draw[dotted, thick, color=black] (3.5,.5) -- (7.5,2.5);
\draw[dotted, thick, color=black] (3.5,1.3) -- (3.5,1.7);
\draw[dotted, thick, color=black] (7.5,1.3) -- (7.5,1.7);
\end{tikzpicture}\\
\caption{A pictorial view of $\bm{\gamma}$ for $\delta_{x_0}$ (thick line $t=0$, dashed $t_1$, dotted $t_2$ with $t_2>t_1>0$)}
\label{fig:pictorialview}
\end{figure}

We will characterize the sub-differential using the following closure property~\cite[Lemma 10.3.8]{AGS}:
\begin{lemma}[Closure of the sub-differential]\label{lemmaclosure}
Let $\phi_h:\mathcal{P}_2(\R) \to (-\infty,+\infty]$ be $\lambda$-geodesically functionals which $\Gamma(d_W)$-converge to $\phi$ as $h \to \infty$. If
\begin{equation*}
\begin{split}
& \bm{\gamma}_h \in \bm{\partial}\phi_h(\mu_h), \qquad \mu_h \to \mu \quad \text{in } \mathcal{P}_2(\R),\qquad \mu \in D(\phi), \\
& \sup_h|\bm{\gamma}_h|_2 < +\infty, \qquad \bm{\gamma}_h \to \bm{\gamma} \quad \text{in } \mathcal{P}_2(\R \times \R),
\end{split}
\end{equation*}
then $\bm{\gamma} \in \bm{\partial}\phi(\mu)$
\end{lemma}
The functional $\W$ has been proven in Proposition \ref{0geodesicconvex} to be 0-geodesically convex, so we can use this lemma with the sequence $\phi_h := \W$ which $\Gamma$-converges to itself. In the following we will use that a measure $\mu_0$ can always be written as $ \mu_0 = \nu + \sum_{i\in I} m_i \delta_{x_i}$ with $\nu(\{x\}) =0$ for every $x \in \R$ for a index set $I$ finite or countable. We define $\alpha_i$ and $\beta_i:=\alpha_i + m_i$ such that $X_{\mu_0}=x_i$ on $(\alpha_i,\beta_i)$. We can now state the following:

\begin{proposition}\label{subdiffp}
Given the functional $\W$ and a measure $\PP \ni \mu_0 = \nu + \sum_{i\in I} m_i \delta_{x_i}$, for some finite or countable $I$ and with $\nu(\{x\})=0$ for every $x \in \R$, then, defining $\Delta_i=\left[ 2\alpha_i-1,2\beta_i -1\right]$, $\mathcal{X}_{\Delta_i}$ the characteristic function of the interval $\Delta_i$ and $k_0(x):=2F_0(x)-1$, the plan
\begin{equation} \label{gamma}
\bm{\gamma}(x,y) =  \sum_i \frac{1}{2}\delta_{x_i} \otimes \mathcal{X}_{\Delta_i} + (i\otimes k_0)_{\#}\nu,
\end{equation}
is the unique element of minimal norm in $\bm{\partial}\W(\mu_0)$.
\proof
Let $\mu_t$ be the Wasserstein gradient flow solution of $\W(\mu)$ with $\sigma=-1$ starting from $\mu_0$; $\mu_t$ is  absolutely continuous for $t > 0$ as we remarked in Section~\ref{sec:particle}. By Theorem~\ref{th:charsubdiff}, for  $t>0$ the extended sub-differential is therefore
\begin{equation}\label{tech44}
\bm{\gamma}_t=(i \otimes k_t)_{\#}\mu_t \quad \text{ with } k_t(x)=2F_t(x)-1,
\end{equation}
writing $F_t$ for the cumulative distribution function of $\mu_t$ as above.
We apply Lemma \ref{lemmaclosure} to the sequences $\mu_t$ and $\gamma_t$ as $t\to 0$.

First note that for every test function $\phi \in C_b(\R^2)$, by the absolute continuity of $\mu_t$ (see Theorem~\ref{thm:dwcontraction}) and the property that $F(X(s))=s$ when the corresponding $\mu$ is absolutely continuous, we have that
\begin{equation*}
\int_{\R^2} \phi(x,y)\,d\bm{\gamma}_t=\int_{\R} \phi(x,2F_t(x)-1)\,d\mu_t= \int_0^1 \phi(X_t(s),2s-1)\,ds.
\end{equation*}
We define $\bm{\gamma}$ through the limit of this expression as $t\to 0$, i.e.
so that, taking the limit for $t \to 0$, which we can calculate explicitly by,
\begin{equation*}
\int_{\R^2} \phi(x,y)d\bm{\gamma} := \lim_{t \to 0}\int_0^1 \phi(X_t(s),2s-1)ds = \int_0^1 \phi(X_0(s),2s-1)ds.
\end{equation*}

First note that $\bm{\gamma}_t \to \bm{\gamma}$ not only narrowly on $\R\times \R$---by construction---but also in the Wasserstein metric on $\R\times \R$, i.e.\ in $\mathcal P_2(\R\times\R)$, since the second moments of $\bm\gamma_t$ also converge (see Figure~\ref{fig:pictorialview} for an illustration).

Next we characterize $\bm\gamma$ in the following way. Set $\Omega_X:= \bigcup_{i\in I}(\alpha_i,\beta_i)$ and then set $Y := [0,1]\setminus \Omega_X$.
Writing
\begin{equation*}
\int_0^1\phi(X_0(s),2s-1)\,ds=\int_{\Omega_X}\phi(X_0(s),2s-1)\,ds+ \int_Y \phi(X_0(s),2s-1)\,ds,
\end{equation*}
the second term can be written as
\begin{equation*}
\int_Y \phi(X_0(s),2s-1)\,ds = \int_{\R  } \phi(x,y)\,d(i \otimes k_0)_{\#}\nu,
\qquad\text{with }k_0(x) = 2F_0(x)-1,
\end{equation*}
and the first as
\begin{align*}
 \int_{\Omega_X}\phi(X_0(s),2s-1)\,ds &=\sum_{i\in I} \int_{F_0(x_i)-m_i}^{F_0(x_i)}\phi(x_i,2s-1)\,ds \\
 &= \sum_i \frac{1}{2} \int_{2F_0(x_i)-2m_i-1}^{2F_0(x_i)-1}\phi(x_i,s')\,ds' = \sum_i \int_{\R} \phi(x,y)\, d\left( \frac{1}{2}\delta_{x_i} \otimes \mathcal{X}_{\Delta_i}\right).
\end{align*}
Therefore
\begin{equation*}
\bm{\gamma} = \sum_i \frac{1}{2}\delta_{x_i} \otimes \mathcal{X}_{\Delta_i} + (i \otimes k_0)_{\#}\nu.
\end{equation*}
We can directly calculate the norm of $\bm{\gamma}$
\begin{equation*}
|\bm{\gamma}|_2^2= \int_0^1 (2s -1)^2ds=\frac{1}{3}\,.
\end{equation*}
It is easy to check by using \eqref{tech44} that $|\bm{\partial} \W|^2(\mu_t)=\tfrac13$. Since the slope is non increasing along solutions we conclude that $|\bm{\partial} \W|^2(\mu_0)\geq\tfrac13$ proving the minimality of $\bm{\gamma}$.
\endproof
\end{proposition}

\begin{remark}
We have that the concept of extended sub-differential is absolutely needed in the repulsive case when dealing with Dirac delta functions. Moreover we just proved that $D(\bm{\partial}\W)$ is the whole $\PP$.
\end{remark}

\section{Discussion}
\label{sec:discussion}

The results of this paper create a connection between two systems that are individually well-studied but are often considered completely separate: the entropy-solution interpretation of conservation laws on one hand and the metric-space gradient flows on the other. For both systems, smooth solutions are unique and reversible in time, but for non-smooth solutions both the uniqueness and the reversibility break down---and these two issues are strongly connected through the concept of information loss.

It is intriguing to see that the gradient-flow concept singles out the same solution as the Oleinik-Kru\v zkov entropy condition, thereby distinguishing the solution from other types with `non-classical' shocks. When dealing with the attractive case, the Rankine-Hugoniot condition uniquely characterizes shocks by our choice of increasing initial conditions. But in the repulsive case non-entropic shocks can occur, and we recover uniqueness with the Oleinik entropy condition. The non-uniqueness for the conservation law is translated into non-uniqueness for \eqref{eq:main}, in a one-to-one correspondence, where persistence vs. spreading of a shock translates into  persistence vs. spreading of a Dirac delta function.

The question naturally arises how the gradient-flow concept embodies the same selection criterion as the entropy condition. Or, to phrase it differently, since both solution concepts lead to non-reversibility in time, how does this non-reversibility arise? One way to see this is the fact that both solution concepts contain a specific {inequality}. In the Burgers equation the inequality is explicitly given in the definition (see~\eqref{Oleinik}); in the gradient-flow solution the inequality lies in the fact that solutions are \emph{curves of steepest descent}. Since forward-time and backward-time evolutions differ by the sign of the functional, the condition of steepest \emph{descent} distinguishes between the two.

In the case of the gradient-flow concept, the non-reversibility seems tightly connected to the Fr\'echet sub-differential and the metric slope at any given time. In the attractive case and for purely atomic initial data, the metric slope has a decreasing jump at each collision time between two particles. This is how the irreversibility shows in the system and it is determined by the element of minimal norm in the Fr\'echet sub-differential. In the repulsive case, the Fr\'echet sub-differential for an initial data with concentrated mass is empty, and thus the system chooses the velocity distribution for the concentrated mass with uniform probability on the admissible velocity range. This is expressed mathematically by the explicit formula of the element of minimal norm in the extended subdifferential given in Proposition~\ref{subdiffp}.

Finally, let us remark that this equivalence is very specific for the attractive/repulsive Newtonian potentials in one dimension, since by integrating the nonlocal equation \eqref{eq:main}, we usually get a nonlocal conservation law except for $W(x)=\pm |x|$.  In other words, the flux of the conservation law can only be expressed as an explicit function of the cumulative distribution function for these two specific cases.

\section*{Acknowledgments}
We thank Giuseppe Savar\'e and Upanshu Sharma for the precious comments and help given during the preparation of this work. A particular thank goes to the whole CASA group from Technische Universiteit Eindhoven for the interesting and helpful discussions. MDF is supported by the FP7-People Marie Curie CIG (Career Integration Grant) Diffusive Partial Differential Equations with Nonlocal Interaction
in Biology and Social Sciences (DifNonLoc), by the `Ramon y Cajal' sub-programme (MICINN-RYC) of the Spanish Ministry of Science and
Innovation, Ref. RYC-2010-06412, and by the by the Ministerio de Ciencia e Innovaci\'on, grant MTM2011-27739-C04-02. JAC was partially supported by the project
MTM2011-27739-C04-02 DGI (Spain), by the 2009-SGR-345 from
AGAUR-Generalitat de Catalunya, by the
Royal Society through a Wolfson Research Merit Award, and by the Engineering and Physical Sciences Research Council grant number EP/K008404/1.
GAB and MAP kindly acknowledge support from the Nederlandse Organisatie voor Wetenschappelijk Onderzoek (NWO) VICI grant 639.033.008.

\medskip

\bibliography{BCDP_arXiv}
\bibliographystyle{plain}

\end{document}